\documentclass[12pt]{amsart}
\usepackage{amscd,amssymb,showkeys}

\usepackage{float}

\theoremstyle{latex 2e}
\newtheorem{thm}[subsection]{Theorem}
\newtheorem{lem}[subsection]{Lemma}
\newtheorem{prop}[subsection]{Proposition}

\newtheorem{defn}[subsection]{Definition}

\theoremstyle{remark}

\numberwithin{equation}{section}

\begin{document}
\title
{Maximum solutions of normalized Ricci flows on 4-manifolds }

\author[F. Fang]{Fuquan Fang}
\thanks{The first author was supported by
NSF Grant 19925104 of China, 973 project of Foundation Science of
China, and the Capital Normal University}
\address{Department of Mathematics, Capital Normal University,
Beijing, P.R.China}
  \email{ffang@nankai.edu.cn}
\author[Y. Zhang]{Yuguang Zhang}
\address{Department of Mathematics, Capital Normal University,
Beijing, P.R.China  }
\author[Z. Zhang]{Zhenlei Zhang}
\address{Nankai Institute of Mathematics,
Weijin Road 94, Tianjin 300071, P.R.China}

\begin{abstract}
We consider maximum solution $g(t)$, $t\in [0, +\infty)$, to the
normalized  Ricci flow. Among other things, we prove that, if $(M,
\omega) $ is a smooth  compact symplectic $4$-manifold such that
$b_2^+(M)>1$ and let $g(t),t\in[0,\infty)$, be a solution to (1.3)
on $M$ whose Ricci curvature satisfies that
$|\text{Ric}(g(t))|\leq 3$ and additionally $\chi(M)=3 \tau
(M)>0$, then  there exists an $m\in \mathbb{N}$, and  a sequence
of points $\{x_{j,k}\in M\}$, $j=1, \cdots, m$, satisfying that,
by passing to a subsequence,
$$(M, g(t_{k}+t),
x_{1,k},\cdots, x_{m,k}) \stackrel{d_{GH}}\longrightarrow ( \coprod
_{j=1}^m N_j , g_{\infty}, x_{1,\infty}, \cdots, , x_{m,\infty}),$$
$t\in [0, \infty)$,  in the $m$-pointed Gromov-Hausdorff sense for
any sequence $t_{k}\longrightarrow \infty$, where  $(N_{j},
g_{\infty})$, $j=1,\cdots, m$, are complete complex hyperbolic
orbifolds of complex  dimension 2 with at most finitely many
isolated orbifold points. Moreover, the convergence is $C^{\infty}$
in the non-singular part of $\coprod _1^m N_{j}$ and
$\text{Vol}_{g_{0}}(M)=\sum_{j=1}^{m}\text{Vol}_{g_{\infty}}(N_{j})$,
where $\chi(M)$ (resp. $\tau(M)$) is the Euler characteristic (resp.
signature) of $M$.
\end{abstract}

\maketitle

\section{Introduction}

Let $(M,g)$ be a compact Riemannian manifold. The Perelman
$\lambda$-functional
\begin{equation}\lambda_{M}(g)=\inf _{f\in C^\infty(M)} \{\mathcal{F}(g,f): \int_{M}e^{-f}
d\text{vol}_{g}=1\}\end{equation} where $\mathcal{F}(g,f)=\int_{M}
(R_g+|\nabla f|^{2})e^{-f} d\text{vol}_{g}$ and $R_g$ is the
scalar curvature of $g$. Note that $\lambda_{M}(g)$ is the lowest
eigenvalue of the operator $-4\triangle+R_{g}$. By [Pe1] the
gradient flow of the Perelman $\lambda$-functional is the
Hamilton's the Ricci-flow evolution equation
\begin{equation}
\frac{\partial}{\partial t}g(t)=-2{\rm Ric}(g(t))
\end{equation}
The normalized Ricci flow equation on an $n$-manifold $M$ reads
\begin{equation}
\frac{\partial}{\partial
t}g(t)=-2Ric(g(t))+\frac{2\overline{R}}{n}g(t)\end{equation}
 where
$Ric$ (resp. $\overline{R}$) denotes the Ricci tensor (resp. the
average scalar curvature $\frac{\int_{M}Rdv}{\int_{M}dv}$). Note
that (1.2) and (1.3) differ only by a change of scale in space and
time, and the volume $\text{Vol}(g(t))$ is constant in $t$. If
$\text{dim}M=n$, $\overline{\lambda}_{M}(g)
=\lambda_{M}(g)\text{Vol}_{g}(M)^{\frac{2}{n}}$ is invariant up to
rescaling the metric. Perelman [Pe1] has proved that
$\overline{\lambda} _M(g(t))$ is non-decreasing along the Ricci
flow $g(t)$ whenever $\overline{\lambda} _M(g(t))\leq 0$. This
leads to the Perelman invariant $\overline{\lambda}_{M}$ by taking
supremum of $\overline{\lambda}_{M}(g)$ in the set of all
Riemannian metrics on $M$.

 By [AIL] the Perelman invariant
$\overline{\lambda }_M$ is equal to the Yamabe invariant whenever
$\overline{\lambda}_{M}\leq 0$, after the earlier estimations (cf.
[An5] [Pe2] [Le4]   [FZ] and [Kot]). In particular, if $(M, g)$ is a
smooth compact oriented $4$-manifold with a  $\rm
Spin^{c}$-structure $\mathfrak{c}$ which is a monopole class (i.e.,
the associated Seiberg-Witten equation possesses an irreducible
solution) so that that $c_{1}^{2}(\mathfrak {c})[M]
>0$, by [FZ] $\overline{\lambda}_{M}\leq - \sqrt{32\pi^{2}
c_{1}^{2}(\mathfrak {c})[M]}.$  Moreover, $g$ is a
K\"{a}hler-Einstein metric of negative scalar curvature if and only
if  $\overline{\lambda}_{M}(g)=- \sqrt{32\pi^{2} c_{1}^{2}(\mathfrak
{c})[M]}$. However, there are plenty of $4$-manifolds where the
Perelman invariant $\overline{\lambda }_M=- \sqrt{32\pi^{2}
c_{1}^{2}(\mathfrak {c})[M]}$ but do not admit any K\"ahler Einstein
metric. It is natural to study $4$-manifolds with these extremal
property. For such a $4$-manifold $M$, to seek for an "optimal"
Riemannian metric on $M$ with respect to the Perelman functional
$\overline \lambda _M: \mathcal M\to \Bbb R$, we want to consider a
{\it maximal solution} $g(t)$ which is a solution of the Ricci flow
(1.3).  We call a longtime solution $g(t)$, $t\in [0, +\infty)$, to
the Ricci flow (1.3) a {\it maximum solution} if $\lim\limits_{t\to
\infty}
 \overline{\lambda}_{M}(g(t))=\overline{\lambda}_{M}$. For a compact $3$-manifold, by
Perelman [Pe2] all solutions of the Ricci flow (1.2) with surgery
exist for longtime and are maximum solutions, provided
$\overline{\lambda}_{M}\leq 0$. In the paper [FZZ] obstructions are
found for the longtime solutions with bounded curvature to (1.3).

In this paper we are going to study the maximum solutions of (1.3)
with bounded Ricci curvatures instead.  To avoid technique
terminology we only state our results for symplectic $4$-manifolds
by using the celebrated work of Taubes [Ta]: if $(M, \omega)$ is a
compact symplectic manifold with $b_2^+(M)>1$ (the dimension of
self-dual harmonic $2$-forms of $M$), the spin$^c$-structure
induced by $\omega$ is a monopole class. Moreover, in this
situation $c_1^2(\mathfrak {c}) [M]=2\chi (M)+3\tau (M)$, where
$\chi(M)$ (resp. $\tau(M)$) is the Euler characteristic (resp.
signature) of $M$.

\vskip 3mm

\begin{thm}  Let $(M, \omega )$ be a smooth  compact
symplectic $4$-manifold satisfying that $b_2^+(M)>1$ and $2\chi
(M)+3\tau (M)>0$. If $g(t),t\in[0,\infty)$, is a solution to (1.3)
such that $|Ric(g(t))|\leq 3$, and
$$\lim\limits_{t\to \infty}\overline{\lambda}_{M}(g(t))=-
\sqrt{32\pi^{2}(2\chi (M)+3\tau (M))},$$ then there exists an $m\in
\mathbb{N}$, and   sequences  of points $\{x_{j,k}\in M\}$, $j=1,
\cdots, m$, satisfying that, by passing to a subsequence,
$$(M, g(t_{k}+t),
x_{1,k},\cdots, x_{m,k}) \stackrel{d_{GH}}\longrightarrow ( \coprod
_{j=1}^m N_j , g_{\infty}, x_{1,\infty}, \cdots, , x_{m,\infty}),$$
 $t\in [0, \infty)$, in the $m$-pointed Gromov-Hausdorff sense for
any sequence $t_{k}\longrightarrow \infty$, where  $(N_{j},
g_{\infty})$, $j=1,\cdots, m$, are complete
 K\"ahler-Einstein orbifolds of complex  dimension 2 with at most
finitely many isolated orbifold points.
 The     scalar curvature (resp. volume) of $g_{\infty}$ is
$$-\text{Vol}_{g_{0}}(M)^{-\frac{1}{2}}\sqrt{32\pi^{2}
(2\chi (M)+3\tau (M))} \ \ \ {\rm (resp. }  \
\text{Vol}_{g_{0}}(M)=\sum_{j=1}^{m}\text{Vol}_{g_{\infty}}(N_{j}))$$
Moreover, the convergence is $C^{\infty}$ in the non-singular part
of $\coprod _1^m N_{j}$.
\end{thm}

We first remark that, if the diameters $\text{diam}_{g(t_{k})}(M)$
possess a uniform upper bound, then $m=1$, and $N_{1}$ is a
compact K\"ahler-Einstein orbifold. Secondly, if the Ricci
curvature bound in the above theorem is replaced by a uniform
bound of sectional curvature, then every $(N_{j}, g_{\infty})$,
$j=1, \cdots ,m$ are complete K\"ahler-Einstein manifolds. By the
same arguments as in [An5][An6], $\coprod_{j=1}^{m}N_{j}$ can
weakly  embed in  $M$, $\coprod _{j=1}^{m}N_{j}
 \subset\subset M$, i.e. for any compact subset $K\subset
 \coprod_{j=1}^{m}N_{j}$, there is a smooth embedding $F_{K}:
 K\longrightarrow M$.  Furthermore, there exists a sufficiently large compact
 subset $K\subset \coprod_{j=1}^{m}N_{j}$ such that $M\backslash K$ admits an
 F-structure of positive rank. This type geometric decomposition
 seems very useful to understand the diffeomorphism type  of
   $4$-manifolds.

\begin{thm} Let $(M, \omega) $ be a smooth  compact
symplectic $4$-manifold such that $b_2^+(M)>1$ and let
$g(t),t\in[0,\infty)$, be a solution to (1.3) such that
$|R(g(t))|\leq 12$. If in addition $\chi(M)=3 \tau (M)>0$, then
$$\lim\limits_{t\to \infty}\overline{\lambda}_{M}(g(t))=-
\sqrt{32\pi^{2}(2\chi (M)+3\tau (M))}$$  Moreover, if
$|Ric(g(t))|\leq 3$, the K\"ahler-Einstein metric $g_{\infty}$ in
Theorem 1.1 is complex hyperbolic.
\end{thm}

To conclude the section we point out that the main result in Theorem
1.1 (resp. Corollary 1.2) holds if the manifold is not symplectic
but a compact oriented $4$-manifold with a monopole class $c_1$
(i.e. with a spin$^c$-structure with non-vanishing Seiberg-Witten
invariant) so that $c_1^2 =  2\chi (M)+3\tau (M)>0$.

 \vskip 8mm
\section{Preliminaries }

\subsection{Monopole class}

Let $(M, g)$  be a   compact oriented  Riemannian  $4$-manifold
with a  $\rm Spin^{c}$ structure $\mathfrak{c}$. Let
$b^{+}_{2}(M)$ denote the dimension of the space of self-dual
harmonic $2$-forms in $M$. Let $S^{\pm}_{\mathfrak{c}}$ denote the
$\rm Spin^{c}$-bundles associated to $\mathfrak{c}$, and let
 $L$ be the determinant line bundle of $\mathfrak{c}$. There is a well-defined
 Dirac operator $$\mathcal{D}_{A}:
  \Gamma(S^{+}_{\mathfrak{c}})\longrightarrow
  \Gamma(S^{-}_{\mathfrak{c}})$$

Let $c:
  \wedge^{*}T^{*}M \longrightarrow {\rm End}(S^{+}_{\mathfrak{c}}\oplus
  S^{-}_{\mathfrak{c}})$ denote the Clifford multiplication on the $\rm{Spin}^c$-bundles,
and, for any $\phi\in \Gamma(S^{\pm}_{\mathfrak{c}})$, let
  $$q(\phi)=\overline{\phi}\otimes\phi-\frac{1}{2}|\phi|^{2}{\rm id}.$$
   The Seiberg-Witten equations read
   \begin{equation}\begin{array}{ccc}\mathcal{D}_{A}\phi=0 \\
   c(F^{+}_{A})=q(\phi)
\end{array}   \end{equation}
 where $A$ is  an  Hermitian
 connection  on $L$, and $F^{+}_{A}$ is the self-dual
 part of the curvature of $A$.

 A solution of (2.1) is called  {\it reducible} if $\phi\equiv 0$;
 otherwise, it is called {\it irreducible}. If $(\phi, A)$ is a resolution of (2.1), one calculates
 easily that  \begin{equation}|F^{+}_{A}|=\frac{1}{2\sqrt{2}}
 |\phi|^{2},  \end{equation}

The Bochner formula reads
\begin{equation} 0=-2\triangle |\phi|^{2}
 +4|\nabla^{A}\phi|^{2}+R_{g}|\phi|^{2}+|\phi|^{4}, \end{equation} where
 $R_{g}$ is the scalar curvature of $g$.

The Seiberg-Witten invariant can be defined by counting the
irreducible solutions of the Seiberg-Witten equations (cf. [Le2]).

\begin{defn}\rm ([K1]) Let $M$ be a smooth compact
  oriented   $4$-manifold. An element
  $\alpha\in H^{2}(M, \mathbb{Z})/$torsion is called a monopole
  class of $M$ if and only if there exists a $\rm Spin^{c}$-structure
  $\mathfrak{c}$ on $M$ with first Chern class $c_{1}\equiv
  \alpha$\rm(mod torsion), \it so that the Seiberg-Witten equations have a solution for every
   Riemannian metric $g$ on $M$.
\end{defn}

By the celebrated work of Taubes [Ta], if  $(M, \omega)$ is a
compact symplectic 4-manifold with $b^{+}_{2}(M)>1$, the canonical
class of $(M, \omega)$ is a monopole class.

\vskip 5mm
\subsection{Kato's
inequality} Let $(M, g)$ be a  Riemannian  $\rm Spin^{c}$-manifold
of dimension $n$, the following Kato inequality is useful.

\begin{prop}\rm(Proposition 2.2 in [BD])  \it Let $\phi$ be a harmonic $\rm Spin^{c}$-spinor on $(M, g)$, i.e. $\mathcal{D}_{A}\phi=0$, where
 $\mathcal{D}_{A}$ is the Dirac operator and $A$ is an  Hermitian  connection on
 the determinant line bundle. Then \begin{equation}|\nabla|\phi||^{2}\leq \frac{n-1}{n}|\nabla^{A}\phi|^{2}\leq |\nabla^{A}\phi|^{2}
 \end{equation} at all points where $\phi$ is non-zero. Moreover,
  $|\nabla|\phi||^{2}= |\nabla^{A}\phi|^{2}$ occurs only if $\nabla^{A}\phi\equiv 0$.
\end{prop}

  Note that  the arguments in
the proof of Proposition 2.2 in [BD] can be used to prove this
proposition without any change, where the same conclusion was
derived for $\rm Spin$-spinor $\phi$.
      For any $\epsilon >0$, let
      $|\phi|_{\epsilon}^{2}=|\phi|^{2}+\epsilon^{2}$. If
      $\phi$ is harmonic, by above proposition, \begin{equation}
      |\nabla|\phi|_{\epsilon}|^{2}\leq \frac{|\phi|}{|\phi|_{\epsilon}}|\nabla|\phi||^{2}
      \leq \frac{n-1}{n}|\nabla^{A}\phi|^{2}\leq |\nabla^{A}\phi|^{2} \end{equation} at points where
      $\phi(p)\neq 0$.  Since $\{p\in M:  \phi(p)\neq 0\} $ is dense in $M$ for harmonic
      $\phi$, we conclude that (2.5) holds  everywhere in $M$.

\vskip 5mm
\subsection{Chern-Gauss-Bonnet formula and  Hirzebruch signature formula} Let $(M,g)$ be a compact
 closed oriented  Riemannian  4-manifold, $\chi(M)$ and $\tau(M)$ are
the Euler number and the  signature of $M$ respectively.   The
Chern-Gauss-Bonnet formula and the Hirzebruch signature theorem say
that \begin{equation}\chi(M)=
\frac{1}{8\pi^{2}}\int_{M}(\frac{R_{g}^{2}}{24}+|W_{g}|^{2}-\frac{1}{2}
|Ric\textordmasculine|^{2})dv_{g}, \ \ \ \ { \rm
and}\end{equation}\begin{equation}
\tau(M)=\frac{1}{12\pi^{2}}\int_{M}(|W^{+}_{g}|^{2}-|W^{-}_{g}|^{2})dv_{g},\end{equation}
where $Ric\textordmasculine =Ric(g)-\frac{R_{g}}{4}g$ is the
Einstein tensor,  $W^{+}_{g}$ and $W^{-}_{g}$ are the self-dual and
anti-self-dual Weyl tensors respectively (cf. [B]). If $g$ is a
K\"ahler-Einstein metric, then
\begin{equation}R_{g}^{2}=24|W^{+}_{g}|^{2},\end{equation} (cf. [Le3]) which will be used in
 the proof of Theorem 1.1.

  By Chern-Gauss-Bonnet formula, one has an $L^{2}$-bound of the curvature operator $Rm(g)$ by the
bounds of Ricci curvature, i.e. if $|Ric(g)|<C$, then
\begin{equation}\int_{M}|Rm(g)|^{2}dv_{g}\leq 8\pi^{2}\chi(M)+C_1 Vol_{g}(M),
\end{equation} where $C$ and $C_1$ are  constants  independent of $(M,g)$.

Let $(N, g)$ be a complete Ricci-flat Einstein 4-manifold. Assume
that
\begin{equation}\int_{N}|Rm(g)|^{2}dv_{g}<\infty, \ \ \ {\rm and}
\ \ \ \rm \text{Vol}_{g}(B_{g}(x,r))\geq Cr^{4},
\end{equation} for all $r> 0$, a point $x\in N$, and a positive
constant $C$.  By Theorem 2.11 of [N],  $(N, g)$ is ALE. (i.e,
Asymptotically Locally Euclidean space) of order 4. It is
well-known that $N$ is asymptotic to the cone on the spherical
space form $S^{3}/ \Gamma$, where $\Gamma \subset SO(4)$ is a
finite group. The Chern-Gauss-Bonnet formula implies that
\begin{equation}\chi(N)=\frac{1}{8
    \pi^{2}}\int_{N}|Rm(g)|^{2}dv_{g}+
    \frac{1}{|\Gamma|}\end{equation} (cf. [N] and [An1]).

\subsection{Curvature estimates  for  4-manifolds}  Now let's recall
a result of [CT], which is important to the  proof of Theorem 1.1.
Let $(M,g)$ be a complete  Riemannian  4-manifold. A subset
$U\subset M$ such that for all $p\in U$,
$\sup\limits_{B_{g}(p,1)}Ric(g)\geq -3$, is called
$\varrho$-collapsed if for all $p\in U$,
$$Vol_{g}(B_{g}(p,1))\leq \varrho.$$
By Theorem 0.1 in [CG], there is a constant $\varepsilon_{4}$ such
that if $U$ is $\varrho$-collapsed with sectional  curvature
$|K_{g}|\leq 1$ and $\varrho \leq \varepsilon_{4}$, then $U$
carries an F-structure of positive rank.

\begin{thm}\rm(Remark 5.11 and Theorem 1.26 in [CT]) \it There exist constants  $\delta >0$, $c>0$ such
that: if $(M, g)$ is a complete  Riemannian  4-manifold with
$|Ric(g)|\leq 3$ and $$\int_{M}|Rm(g)|^{2}dv_{g}\leq C,$$ and if
$E\subset M$ is a bounded open subset such that $T_{1}(E)=\{ x\in
M: \text{dist}(x, E)\le 1\}$  is $\varepsilon_{4}$-collapsed with
$$
 \int_{B_{g}(x,1)}|Rm(g)|^{2}dv_{g}\leq \delta  \ \ \ \ ({ for \ \ \ all} \ \ \ \ T_{1}(E)), $$ then
 $$\int_{E}|Rm(g)|^{2}dv_{g}\leq cVol_{g}(A_{0,1}(E)),$$ where $A_{0,1}(E)=T_{1}(E)\backslash E$.  \end{thm}

\vskip 3mm

\section{ The limiting  behavior of Ricci flow }
In this section we study the limiting  behavior of Ricci-flow with
bounded Ricci curvatures on $4$-manifolds. We will assume in this
section that $M$ is a smooth closed oriented $ 4$-manifold with
$\overline{\lambda}_{M}<0$, and  $g(t)$, $t\in [0, +\infty)$, is a
longtime solution of the  normalized Ricci flow (1.3) with bounded
Ricci-curvature. By normalization we may assume that
$|Ric(g(t))|\leq 3$. By (2.9) there is a constant $C$ independent
of $t$  such that
$$\int_{M}|Rm(g(t))|^{2}dv_{g(t)}\leq C.$$

Let us denoted by $V$ the volume
$\text{Vol}_{g(0)}(M)=\text{Vol}_{g(t)}(M)$, and
$\breve{R}(g(t))=\min\limits_{x\in M} R(g(t))(x)$ the minimum of
the scalar curvature of  $g(t)$. It is easy to see that
$\breve{R}(g(t))\leq \overline{\lambda}_{M}V^{-\frac{1}{2}} < 0$.

 \begin{lem}

 (3.1.1) $\lim\limits_{t\rightarrow\infty}\lambda_{M}(g(t))=\lim\limits_{t\to \infty}\overline{R}(g(t))=
 \lim\limits_{t\to \infty}\breve{R}(g(t))=\overline{R}_{\infty}$\newline
 (3.1.2) $ \lim\limits_{t\to \infty}
 \int_{M}|R(g(t))-\overline{R}(g(t))|dv_{g(t)}=0,$\newline
 (3.1.3) $ \lim\limits_{t\to
\infty} \int_{M}|Ric\textordmasculine (g(t))|^{2}dv_{g(t)}=0.$
  \end{lem}
 \begin{proof} By Perelman [Pe1] $\lambda _M(g(t))$ is a
 non-decreasing function on $t$, therefore the limit $\lim\limits_{t\to \infty}\lambda _M(g(t))$
 exists since $\overline \lambda _M<0$. Now let us denote by
 $\overline{R}_{\infty}$ the limit
$\lim\limits_{t\rightarrow\infty}\lambda_{M}(g(t))$. Note that
$\overline{R}_{\infty} \leq \overline{\lambda}_{M}
V^{-\frac{1}{2}}<0$. To prove (3.1.1), we first prove that both
$\lim\limits_{t\rightarrow\infty}\overline{R}(g(t))$ and
$\lim\limits_{t\rightarrow\infty}\breve{R}(g(t))$ exist and take
values $\overline{R}_{\infty}$. By the same arguments as in the
proof of Proposition 2.6 and Lemma 2.7 of [FZZ] we get that
$$\lim\limits_{t\rightarrow\infty}\overline{R}(g(t))-\breve{R}(g(t))=0.$$
Observe that $\overline{R}(g(t))\geq \lambda_{M}(g(t)) \geq
\breve{R}(g(t))$ (cf. [KL] (92.3)). Therefore
$\lim\limits_{t\rightarrow\infty}\overline{R}(g(t))=\overline{R}_{\infty}=\lim\limits_{t\longrightarrow
\infty}\breve{R}(g(t))$. This proves (3.1.1).

Note that
 \begin{eqnarray*}\int_{M}|R(g(t))-\overline{R}(g(t))|dv_{g(t)}& \leq & \int_{M}
 (R(g(t))-\breve{R}(g(t)))dv_{g(t)}+\int_{M}(\overline{R}(g(t))-\breve{R}(g(t)))dv_{g(t)}\\ & = & 2(\overline{R}(g(t))
 -\breve{R}(g(t)))V\end{eqnarray*}
(3.1.2) follows from (3.1.1).

 By Lemma 3.1 in [FZZ], $$\int_{0}^{\infty}\int_{M}|Ric\textordmasculine (g(t))|^{2}dv_{g(t)}dt<\infty,$$
  and, by Lemma 1 in [Ye], we have $$\frac{d}{dt}\int_{M}|Ric\textordmasculine (g(t))|^{2}dv_{g(t)}
 \leq -2\int_{M}|\nabla Ric\textordmasculine (g(t))|^{2}dv_{g(t)} +4 \int_{M}|Rm||Ric\textordmasculine (g(t))|^{2}dv_{g(t)}
 <D,$$ where $D$ is a constant independent of $t$. By the same argument
 as in the proof of Proposition 2.6 in [FZZ] (3.1.3) follows. \end{proof}

The following is the main result of this section, which is an
analogy of Theorem 10.5 in [CT], where the same conclusion was
derived for closed oriented  Einstein $4$-manifolds with the same
negative Einstein constant. The key point in our case is to use
Lemma 3.1 to get non-collapsing balls and  to prove
 the limiting metric is an Einstein metric (cf. Lemma 3.3 and Lemma 3.4 below).

    \begin{prop}  Let $M$ be a smooth closed
oriented $ 4$-manifold with $\overline{\lambda}_{M}<0$. If
$g(t),t\in[0,\infty)$ is a solution to (1.3)  such that
$|Ric(g(t))|\leq 3$, and $\{t_{k}\}$ is a sequence of times tends
to infinity such that
$$\text{diam}_{g_{k}}(M)\longrightarrow \infty,$$ when $k\longrightarrow \infty$, where
$g_{k}=g(t_{k})$, then there exists an $m\in \mathbb{N}$, and
sequences  of points $\{x_{j,k}\in M\}$,  $j=1, \cdots, m$,
satisfying that, by passing to a subsequence,
$$(M, g_{k},x_{1,k},\cdots, x_{m,k}) \stackrel{d_{GH}}\longrightarrow (
\coprod _{j=1}^m N_j , g_{\infty}, x_{1,\infty}, \cdots, ,
x_{m,\infty})$$ in the $m$-pointed Gromov-Hausdorff sense for
${k}\to \infty$, where  $(N_{j}, g_{\infty})$ $j=1,\cdots, m$ are
complete Einstein $4$-orbifolds  with at most finitely many
isolated orbifold points $\{q_{i}\}$.
 The     scalar curvature (resp. volume) of $g_{\infty}$ is
$$\overline{R}_{\infty}=\lim\limits_{t\longrightarrow \infty}\lambda_{M}(g(t)), \ \ \ \ \ {(resp.}  \ \ \ \
V=\text{Vol}_{g_{0}}(M)=\sum_{j=1}^{m}\text{Vol}_{g_{\infty}}(N_{j})).$$
Furthermore, in the regular part of $N_{j}$, $\{g_{k}\}$ converges
to $g_{\infty}$ in  both $L^{2,p}$ (resp. $C^{1,\alpha}$) sense
for all $p<\infty$ (resp. $\alpha <1$).
\end{prop}

We divide the proof of Proposition 3.2 into several useful lemmas.

A key result in the paper [CT] shows that, for any compact
oriented Einstein 4-manifold $(X, g)$ with Einstein constant $-3$,
there exists a constant $C$ depending only on the Euler number of
$X$, and a point $x\in X$ such that $\text{Vol}_{g}  (B_{g}(x,
1))\geq C \text{Vol}_{g} (X)$ (cf. Theorem 0.14 [CT]).
Cheeger-Tian remarked that the same result continues to hold for
$4$-manifolds which are sufficiently negatively Ricci pinched. The
following lemmas is an  analogy of the result for the metric $g_k$
in Proposition 3.2.

 \begin{lem} There exists a constant $v>0$, and a sequence $\{x_{k}\}\subset M$ such that
 $$\text{Vol}_{g_{k}}(B_{g_{k}}(x_{k},1))\geq  v.
 $$
  \end{lem}
    \begin{proof} Let $\varepsilon_{4} >0$ be the critical constant of Cheeger-Tian (cf. $\S$1 [CT]), i.e.,
 if $X$ is a Riemannian $4$-manifold which is $\varepsilon_{4}$-collapsed with
locally bounded curvature, then $X$ carries an  F-structure of
positive rank.  We may assume that, for all $x\in M$ and $g_k$,
$\text{Vol}_{g_{k}}(B_{g_{k}}(x,1))<\varepsilon_{4}$.    By a
standard covering argument, for any $k$, there exist finitely many
points $q_{1}, \cdots, q_{l}$ such that
$E=M\backslash\bigcup_{i=1}^{l} B_{g_k}({q_{i}}, 1)$ satisfies the
hypothesis of Theorem 2.7. Moreover, $l\leq C\delta^{-1}$ where
$C$ and $\delta$ are the constants in Theorem 2.7. Therefore, by
Theorem 2.7 we conclude that, there is a constant $C_1$
independent of $k$ such that
\begin{equation}
\int_{E}|R(g_{k})|^{2}dv_{k}\leq
6\int_{E}|Rm(g_{k})|^{2}dv_{k}\leq
 C_1
 \sum\limits_{i=1}^{l}\text{Vol}_{g_{k}}(B_{g_{k}}(q_{i},1)).
 \end{equation}
On the other hand, by Lemma (3.1.2)  $$|
\int_{E}(R(g_{k})^{2}-\overline{R}(g_{k})^{2})dv_{k}|\leq  24
\int_{E}|R(g_{k})-\overline{R}(g_{k})|dv_{k}\stackrel{k\to\infty}\longrightarrow
0.$$ Therefore
\begin{equation}
\end{equation}
 \begin{eqnarray*}
\frac 12\overline R_\infty^2\text{Vol}_{g_k}(E) -\int
_ER(g_k)^2dv_k
& \leq & \overline R(g_k)^2\text{Vol}_{g_k}(E) -\int _ER(g_k)^2dv_k\\
& = &
\int_{E}(\overline{R}(g_{k})^{2}-{R}(g_{k})^{2})dv_{k}\\
& \leq & \frac 14 \overline R_\infty^2V\end{eqnarray*}
 for sufficiently large $k$ since $\overline
R_\infty\le \overline \lambda _M V^{-\frac12}<0$. By inserting
(3.1) we get that
 \begin{eqnarray*}\frac 12 \overline R_{\infty}^{2}(V-
 \sum\limits_{i=1}^{l}Vol_{g_{k}}(B_{g_{k}}(q_{i},1))-\frac 14 \overline{R}_{\infty}^{2}{V}& \leq
& \frac 12 \overline{R}_{\infty}^{2} \text{Vol}_{g_k}(E)-\frac 14 \overline{R}_{\infty}^{2}{V}\\
& \leq & C_1
 \sum\limits_{i=1}^{l}\text{Vol}_{g_{k}}(B_{g_{k}}(q_{i},1)), \end{eqnarray*}
  and $$V\leq
 C_2\sum\limits_{i=1}^{l}\text{Vol}_{g_{k}}(B_{g_{k}}(q_{i},1)),$$ where  $C_2$ is a constant independent
 of $k$. Therefore, there is at least a ball among the $l$ balls whose volume is at least $\frac V {C_2l}$.
 The desired result follows. \end{proof}

Assuming that $\text{diam}_{g_{k}}(M)\to \infty$ for $k\to
\infty$, by using the technique  developed in [An3], the analogue
of Theorem 3.3 in [An2] holds (cf. Theorem 2.3 in [An4]), i.e.
there exist a sequence
      of points $\{x_{k}\}\subset M$ such that, by  passing to a subsequence,
      $$\{(M, g_{k}, x_{k})\}\stackrel{d_{GH}} \longrightarrow (N_{\infty}, g_{\infty},
      x_{\infty})$$
where $N_{\infty}$ is a $4$-orbifold with only isolated orbifold
points $\{q_{i}\}$, $g_{\infty}$ is a complete  $C^{0}$ orbifold
metric, and $g_{\infty}$ is  a $C^{1,\alpha}\cap L^{2,p}$
Riemannian   metric on the regular  part of $N_{\infty}$, for all
$p<\infty$ and $\alpha <1$. Furthermore, $\{g_{k}\}$ converges to
$g_{\infty}$ in the $L^{2,p}$ (resp. $C^{1,\alpha}$) sense on the
regular  part of $N_{\infty}$, i.e. for any $r\gg 1$ and $k$,
      there is a smooth embedding  $F_{k,r}: B_{g_{\infty}}(x_{\infty}, r)\backslash \bigcup_{i}B_{g_{\infty}}(q_{i},
      r^{-1})\subset N_{\infty}\to M$ such that, by passing to  a
      subsequence, $F_{k,r}^{*}g_{k}$ converge to
      $g_{\infty}$ in both       $L^{2,p}$ and $C^{1,\alpha}$ senses.

\begin{lem} $g_{\infty}$ is an Einstein orbifold  metric with scalar curvature $\overline{R}_{\infty}$.
\end{lem}

\begin{proof}
We first prove  that $g_{\infty}$ is
      an Einstein  metric with scalar curvature $\overline{R}_{\infty}$ on the regular  part of
      $N_{\infty}$. Since $F_{k,r}^{*}g_{k}$ converge to
      $g_{\infty}$ in the $L^{2,p}$(resp. $C^{1,\alpha}$)  sense  on $B_{g_{\infty}}(p_{\infty}, r)\backslash \bigcup_{i} B_{g_{\infty}}(q_{i},
      r^{-1})$, for any $r$, by Lemma 3.1,  we obtain that $$0\leq \int_{B_{g_{\infty}}(p_{\infty}, r)\backslash
\bigcup_{i}B_{g_{\infty}}(q_{i},
      r^{-1})}|Ric\textordmasculine
(g_{\infty})|^{2}dv_{\infty}\leq \lim\limits_{k\longrightarrow
\infty} \int_{M}|Ric\textordmasculine (g_{k})|^{2}dv_{k}=0,$$
$$0\leq \int_{B_{g_{\infty}}(p_{\infty}, r)\backslash
\bigcup_{i}B_{g_{\infty}}(q_{i},
      r^{-1})}| R(g_{\infty})-\overline{R}_{\infty}|dv_{\infty}\leq \lim\limits_{k\longrightarrow
\infty} \int_{M}|R(g_{k})-\overline{R}(g_{k})|dv_{k}=0.$$
Therefore $g_{\infty}$ is a $C^{1,\alpha}$ Riemannian metric on
$B_{g_{\infty}}(p_{\infty}, r)\backslash
\bigcup_{i}B_{g_{\infty}}(q_{i}, r^{-1})$ which satisfies the
Einstein equation in the weak sense. By elliptic regularity
theory, $g_{\infty}$ is a smooth Einstein metric with scalar
curvature $\overline{R}_{\infty}$.

Since  $g_{\infty}$ is a $C^{0}$-orbifold metric, i.e.
    for any orbifold point $q_{i}\in N_\infty $, there is a neighborhood
    $U_{i}\cong
   B(0,r) /\Gamma$ of $q_{i}$ such that   $\widetilde{g}_{\infty}$ is
    a $C^{0}$-Riemannian metric on $B(0,r)\subset \mathbb{R}^{4}$
    where  $\Gamma\subset SO(4)$ is a finite subgroup acting freely on $S^3$, and $\widetilde{g}_{\infty}|_{B(0,r)\backslash \{0\}}$ is
    the pull-back metric of $g_{\infty}$. Note that
    $\widetilde{g}_{\infty}$ is a smooth Einstein metric on
    $B(0,r)\backslash \{0\}$ satisfying that
    $\int_{B(0,r)}|Rm(\widetilde{g}_{\infty})|^{2}dv_{\widetilde{g}_{\infty}}<C<\infty$.
    By the arguments as  in [An1] and [Ti], $\widetilde{g}_{\infty}$ is a
    $C^{\infty}$ Einstein metric on
    $B(0,r)$ (cf. the proof of Theorem C in [An1], and Section 4 in [Ti]).
     Hence $g_{\infty}$ is an Einstein orbifold  metric.
  \end{proof}

By the discussion before Lemma 3.4 we may choose $\ell$ sequences
of points $\{x_{j,k}\}\subset M$, $j=1, \cdots, \ell$, such that
$\text{dist}_{g_{k}}
 (x_{i,k},x_{j,k})\stackrel{k\to \infty}\longrightarrow \infty$ for any $i\neq j$, and
\begin{equation}
\{(M, g_{k}, x_{1,k},
 \cdots, x_{\ell,k}) \}\stackrel{d_{GH}} \longrightarrow (\coprod_{j=1}^{\ell }N_{j}, g_{\infty}, x_{1,\infty},
 \cdots,  x_{\ell,\infty})
 \end{equation}
 where $(N_{j}, g_{\infty}, x_{j,\infty})$, $j=1, \cdots, \ell$ are complete Einstein  $4$-orbifolds with only isolated
 singular points and  scalar curvatures $\overline{R}_{\infty}$.
     Furthermore, $\{g_{k}\}$ converges to $g_{ \infty}$ in  both
      $L^{2,p}$ (resp. $C^{1,\alpha}$) sense on the    regular
      parts  of $N_{j}$, $j=1, \cdots, \ell$. Note that
 \begin{equation}
 V\geq \sum_{i=1}^{\ell }\text{Vol}_{g_{\infty}}(N_{j}).
 \end{equation}

       \begin{lem} The number of orbifold points of
  $ \coprod\limits_{j=1}^{\ell }N_{j}$ is less than a constant depending only on the Euler
  characteristic $\chi (M)$.
    \end{lem}

    \begin{proof} For each  orbifold point $q\in N_{j}$,  there exist a sequence  $\{q_{k}\}\subset M$, and two
constants $r\gg r_{1}>  0$ such that:

(3.5.1)  $q\in    B_{g_{\infty}}(x_{j,\infty}, r)$;

(3.5.2) $B_{g_{\infty}}(q, r_{1})\backslash B_{g_{\infty}}(q,
    \sigma)$ lies in the regular part of $B_{g_{\infty}}(x_{j,\infty},
    r)$ for any $\sigma  <r_{1}$;

(3.5.3) $(B_{g_{k}}(q_{k}, r_{1})\backslash B_{g_{k}}(q_{k},
    \sigma), g_{k})\stackrel{C^{1,\alpha}}\longrightarrow (B_{g_{\infty}}(q, r_{1})\backslash B_{g_{\infty}}(q,
    \sigma), g_\infty)$.

By the definition of harmonic radius (cf. [An3]),  the harmonic
radii of all points in $B_{g_{k}}(q_{k}, r_{1})\backslash
B_{g_{k}}(q_{k},    \sigma)$ have a uniform lower bound, saying
$\mu >0$,  a constant  depending on $\sigma$ but independent of
$k$.

Clearly, there is a positive constant $v_0$ (e.g., $\frac 12
\text{Vol}_{g_\infty} (B_{g_{\infty}}(x_{j,\infty}, r))$) such that
\\ $\text{Vol}_{g_{k}}(B_{g_{k}}(x_{j,k}, r))\geq v_0$. Note that
  the Sobolev constants $C_{S,k}$ of $B_{g_{k}}(x_{j,k}, r)$
are bounded from below by a constant depending only on $v_0, r$ (cf.
[An2] and [Cr]). Therefore, by [An2] again we get that
$\text{Vol}_{g_{k}}(B_{g_{k}}(q_{k}, s))\geq C s^{4}$ for any $s\ll
1$, where $C$ is independent of $k$.

Let us denote by $ r_{h,k}$ the infimum of the harmonic radii of
$g_k$ in the ball $B_{g_{k}}(q_{k},  r_{1})$.  Note that $ r_{h,k}
\stackrel {k\to \infty}\longrightarrow 0$ since $q$ is a orbifold
point (cf. [An3]). Therefore, there is a point $\bar q_{k}\in
B_{g_{k}}(q_{k},
    \sigma)$ so that $r_{h}(\bar q_{k})= r_{h,k}$ for sufficiently large $k$.

Consider the normalized balls $(B_{g_{k}}(q_{k},
    r_{1}), r_{h,k}^{-2}g_{k})$, which have harmonic radii at least $1$.
By passing to a subsequence if necessary,
$$(B_{g_{k}}(q_{k},
    r_{1}),
    r_{h,k}^{-2}g_{k},\bar q_{k})\stackrel{C^{1,\alpha}}\longrightarrow
(W,    \bar g_{\infty}, \bar q)$$ where $(W,\bar g_{\infty}) $ is
a complete Ricci-flat $4$-manifold satisfying that
\begin{equation}
\text{Vol}_{\bar g_{\infty}}(B_{\bar g_{\infty}
    }(\bar q, r))\geq C r^{4}
\end{equation}
for any $r>0$.
      It is obvious that
     $$\int_{W}|Rm(\bar g_{\infty})|^{2}dv_{\bar g_{\infty}}\leq
     \liminf\limits_{k\longrightarrow
     \infty}\int_{M}|Rm(g_{k})|^{2}dv_{k}\leq C.$$ Therefore
 $(W,\bar  g_{\infty})$ is an Asymptotically Locally Euclidean space (cf.
    Theorem 2.11 in [N] or [An1]), which is asymptotic to a cone of
 $S^{3}/ \Gamma$ where  $\Gamma\subset SO(4)$ is a finite group acting freely on $S^3$.
By the Chern-Gauss-Bonnet formula
\begin{equation}
\chi(W)=\frac{1}{8
    \pi^{2}}\int_{W}|Rm(\bar g_{\infty})|^{2}dv_{\bar g_{\infty}}+
    \frac{1}{|\Gamma|}.
\end{equation}
By [An1] $W$ is isometric to $\Bbb R^4$, provided $|\Gamma |=1$.
Since the harmonic radius of $\bar g_\infty$ at $\bar q$ is $1$,
hence $\bar g_\infty$ can not be the Euclidean metric. Hence
$|\Gamma |\geq 2$. It is easy to verify that $\chi (W)\ge 1$. By
(3.6) we get that
    $$\int_{W}|Rm(\bar g_{\infty})|^{2}dv_{\bar g_{\infty}}\geq 4 \pi^{2}.$$
This proves that every orbifold point contributes to
$\liminf\limits_{k\longrightarrow
     \infty}\int_{M}|Rm(g_{k})|^{2}dv_{k}$ at least $4\pi ^2$. By
the rescaling invariance of the integral we conclude that the
number of orbifold points $\beta \leq \frac  C{4 \pi^{2}}$.
\end{proof}

The following lemma is an analogue of a result in Cheeger-Tian
[CT].

  \begin{lem}    $\ell < \chi(M)+ \beta +1$, where $\beta:=\#\{\text{number of
  orbifold points in Lemma 3.5}\}$.   \end{lem}

\begin{proof} Suppose not, i.e, $\ell \geq \chi(M)+ \beta +1 $, by definition there are
at least $\chi(M)+1$ components of  $\coprod _{1}^\ell N_{j}$
which are smooth complete non-compact Einstein $4$-manifolds of
finite volume, for simplicity saying $N_1, \cdots, N_s$, where
$s\geq \chi (M)+1$. By Theorem 4.5 in [CT], for each $1\le j\le
s$,
$$\int_{N_{j}}|Rm(g_{\infty})|^{2}dv_{g_{\infty}}\geq  8\pi^{2}.$$
Since $(M, g_{k}, x_{k,j})\stackrel{L^{2,p}}\longrightarrow (N_{j},
g_{\infty}, x_{\infty,j})$, by Chern-Gauss-Bonnet formula and
(3.1.3) in Lemma 3.1  we get that
$$8\pi^{2} \chi(M)=\lim\limits_{k\longrightarrow
\infty}\int_{M}|Rm(g_{k})|^{2}dv_{g_{k}}\geq \sum
\int_{N_{j}}|Rm(g_{\infty})|^{2}dv_{g_{\infty}}\geq 8\pi^{2}(\chi
(M)+1).$$ A contradiction.
\end{proof}

Let $m$ denote the maximal value of all possible choice of the base
point sequences in (3.3), which has a upper bound by Lemma 3.6.

 \begin{lem}Let $M_{k,r}=M\backslash
\bigcup_{j=1}^{m}B_{g_{k}}(x_{j,k},r)$. For sufficiently large
$r$, there is a constant $C$ independent of $r$ such that
\begin{equation}\lim\limits_{k\to
\infty}\text{Vol}_{g_{k}}(M_{k,r})\leq
 C \sum_{j=1}^{m}\text{Vol}_{g_{\infty}}
 (N_{j}\backslash B_{g_{\infty}}(x_{j,\infty},\frac{r}{2})),\end{equation}
\begin{equation}
 \sum_{j=1}^{m}\text{Vol}_{g_{\infty}}(N_{j})=V.\end{equation}
\end{lem}

  \begin{proof}We may choose $r\gg 1$ such that, for any
  $y \in \coprod _{j=1}^{m}(N_{j}\backslash
   B_{g_{\infty}}(x_{j,\infty},r-1))$, $\rm Vol_{g_{\infty}}(B_{g_{\infty}}(y, 1))\leq
   \frac{1}{2}\varepsilon_{4}$, where $\varepsilon_{4} >0$ is  the critical constant of
Cheeger-Tian (cf. proof of Lemma 3.3 or $\S$1 [CT] ).

Now we claim that there is a constant $k_{0}\gg 1$ such that, for
any $k> k_{0}$ and  any $x\in M_{k,r}$,
$\text{Vol}_{g_{k}}(B_{g_{k}}(x,1))\leq \varepsilon_{4}$.

    If it is false,   without loss of generality we may assume
a sequence of points  $\{y_{k}\}\subset M_{k,r}$ such that
\begin{equation} \rm Vol_{ g_{k}}(B_{g_{k}}(y_{k},1))>
\varepsilon_{4}
\end{equation}
Observe that the distance $\rm dist_{g_{k}}(y_{k},x_{j, k})\to
\infty$ as $k\to \infty$ for all $1\le j\le m$. Otherwise,
assuming $\rm dist_{g_{k}}(y_{k},x_{j, k})< {\rho }$ for some $j$
and $\rho >0 $, we get that $F_{j,k,\rho }^{-1}(y_{k})\to
y_{\infty}\in B_{g_{\infty}}(x_{j,\infty},\rho )\backslash
B_{g_{\infty}}(x_{j,\infty}, r-1)$, and so
\begin{equation}
\rm Vol_{g_{k}}(B_{g_{k}}(y_{k},1))\to \rm
Vol_{g_{\infty}}(B_{g_{\infty}}(y_{\infty},1))\le \frac 12
\varepsilon_{4}
\end{equation}
when $k\to \infty$, since $F_{j,k, \rho}^{*}g_{k} $
$C^{1,\alpha}$-converges  to $g_{j,\infty}$, where
\begin{equation}
F_{j,k, \rho}: B_{g_{\infty}}(x_{j, \infty}, \rho )\backslash
\bigcup _{i}B_{g_{\infty}}(q_{i},
      \rho ^{-1})\subset N_{\infty}\to M
\end{equation}
is a smooth embedding so that $F_{j, k,\rho}^{*}g_{k}$ converges to
$g_{\infty}$ in the $C^{1,\alpha}$-sense (cf. the discussion before
Lemma 3.4). A  contradiction to (3.9).

Note that $(M,g_{k}, y_{k})\stackrel{d_{GH}}\longrightarrow
(N_{\infty},g_{\infty},y_{\infty})$ where $N_\infty$ is a complete
$4$-orbifold different from each of $N_{j}$, $1\le j \le m$. This
violates the choice of maximality of $m$. Hence we have proved the
claim.

By a standard covering argument, for any $k$, there exist finitely
many points $z_{1,k}, \cdots, z_{I,k}$ such
 that  $E_{k,r}=M_{k,r}\backslash\bigcup_{i=1}^{I} B_{z_{i,k}}(1)$ satisfies the hypothesis of Theorem 2.7,
 where $I$ is independent of $k$.  By Theorem 2.7, there is a constant $C$ independent of $k$
 such that
 $$\int_{E_{k,r}}|R(g_{k})|^{2}dv_{k}\leq 6 \int_{E_{k,r}}|Rm(g_{k})|^{2}dv_{k}\leq C
 (\sum\limits_{i=1}^{I}\text{Vol}_{g_{k}}(B_{g_{k}}(z_{i,k},1))+\text{Vol}_{g_{k}}(A_{0,1}(M_{k,r}))).$$
By Lemma 3.1, for $k\gg 1$, we have
\begin{equation}
\int_{E_{k,r}}|R(g_{k})-\overline{R}(g_{k})|dv_{k}<\int_{M}|R(g_{k})-\overline{R}(g_{k})|dv_{k}\longrightarrow
0. \end{equation}
By (3.2) we get
 \begin{eqnarray*}
\frac 12\overline R_\infty^2\text{Vol}_{g_k}({E_{k,r}}) -\int
_{E_{k,r}}R(g_k)^2dv_k & \leq &
\int_{E_{k,r}}(R(g_{k})^{2}-\overline{R}(g_{k})^{2})dv_{k}\\
& \leq  & 24 \int_{E_{k,r}}|R(g_{k})-\overline{R}(g_{k})|dv_{k},
\end{eqnarray*}
Since $\rm Vol_{g_{k}}(E_{k,r})\geq \rm Vol_{g_{k}}(M_{k,r})-
\sum\limits_{i=1}^{I}\rm Vol_{g_{k}}(B_{g_{k}}(z_{i,k},1))$, by
the above together we get immediately that \begin{equation}
\rm
Vol_{g_{k}}(M_{k,r})\leq C
 (\sum\limits_{i=1}^{I}\rm Vol_{g_{k}}(B_{g_{k}}(z_{i,k},1))+\rm Vol_{g_{k}}(A_{0,1}(M_{k,r})))
 +24 \int_{E_{k,r}}|R(g_{k})-\overline{R}(g_{k})|dv_{k}.
\end{equation}
If $\rm dist_{g_{k}}(z_{i,k}, x_{j,k})\to \infty$ for all
 $1\le j\le m$, by the same argument as above we get that  $$\rm Vol_{g_{k}}(B_{g_{k}}(z_{i,k},1))\to 0$$ when
 $k\to \infty$. Otherwise, there exists a  subsequence $k_s\to \infty$ and an index $j$ such that
 $$\rm dist_{g_{k_s}}(z_{i,k_s}, x_{j,k_s})<
 \rho$$ for some constant $\rho$. In both cases, we obtain
 $$\rm \limsup\limits_{k\to\infty}Vol_{g_{k}}(B_{g_{k}}(z_{i,k_s},1))\leq \sum_{j=1}^{m}Vol_{g_{\infty}}
 (N_{j}\backslash B_{g_{\infty}}(x_{j,\infty},\frac{r}{2}))$$
 for $r\gg\rho$. Therefore, by (3.12) and (3.13) we conclude
 immediately (3.7).

 By (3.7) it follows that $\rm  \lim\limits_{k, r\to\infty}Vol_{g_{k}}(M_{k,r})\to
 0$. Hence (3.8) follows.
  \end{proof}

By now Proposition 3.2 follows by the above lemmas.

\vskip 4mm

\section{Smooth convergence on the regular part}

The main result of this section is the following:

\begin{prop}
Let $M$ be a closed $4$-manifold satisfying that
$\bar{\lambda}_{M}<0$ and let $g(t),t\in[0,\infty),$ be a solution
to the normalized Ricci flow equation (1.3) on $M$ with uniformly
bounded Ricci curvature. If $(M,g(t_{k}),p_{k})\stackrel
{d_{GH}}\longrightarrow (N_{\infty},g_{\infty},p_{\infty})$, where
$t_{k}\to \infty$ and $N_\infty$ is a $4$-dimensional orbifold,
and $g(t_k) \stackrel {C^{1,\alpha}}\longrightarrow g_\infty$ on
the regular part $\mathcal{R}$ of $N_\infty$ (the compliment of
the orbifold points), then, by passing to a subsequence, for all
$t\in [0, \infty)$, $(M,g(t_{k}+t),p_{k})   \stackrel
{d_{GH}}\longrightarrow (N_{\infty},g_{\infty}(t),p_{\infty})$,
where $g_{\infty}(t)$ is a family of smooth metrics on
$\mathcal{R}$ solving the normalized Ricci flow equation on
$\mathcal{R}$ with $g_{\infty}(0)=g_{\infty}$. Moreover, the
convergence is smooth on $\mathcal{R}\times [0,\infty)$.
\end{prop}

In [Se] the convergence of K\"ahler-Ricci flow on compact K\"ahler
manifolds with bounded Ricci curvature was studied. It seems that
the arguments in [Se] could be applied to prove Proposition 4.1,
but the authors can not follow completely her line. Therefore, we
give a quite different approach, where we first give a curvature
estimate of the Ricci flow similar to Perelman's pseudolocality
theorem. Using this curvature estimation we prove the limit Ricci
flow exists on $\mathcal{R}\times [0, \infty)$.  Finally, we prove
that $\mathcal R$ is exactly the regular part of every subsequence
limit of $(M,g(t_{k}+t),p_{k})$, for all $t\in [0, \infty)$. It
deserves to point out that our approach works only in dimension
$4$.

We now give a curvature estimate for the Ricci flow which is an
analogy of Perelman's pseudolocality theorem (cf. [Pe1] Thm.
10.1). The difference is that here we use the hypothesis of local
almost Euclidean volume growth, instead of the almost Euclidean
isoperimetric estimate. The proof is much easier than that of
Perelman's pseudolocality theorem.

\begin{thm}
There exist universal constants $\delta_{0},\epsilon_{0}>0$ with
the following property. Let
$g(t),t\in[0,(\epsilon_{P}r_{0})^{2}],$ be a solution to the Ricci
flow equation (1.2) on a closed $n$-manifold $M$ and $x_{0}\in M$
be a point. If the scalar curvature
$$R(x,t)\geq-r_{0}^{-2}\mbox{ whenever }{\rm{dist}}_{g(t)}(x_{0},x)\leq
r_{0},$$ and the volume
$${\rm{Vol}}_{g(t)}(B_{g(t)}(x,r))\geq(1-\delta_{0}){\rm{Vol}}(B(r))\mbox{ for all }B_{g(t)}(x,r)\subset B_{g(t)}(x_{0},r_{0}),$$ where
$B(r)$ denotes a ball of radius $r$ in the $n$-Euclidean space and
${\rm{Vol}}(B(r))$ denotes its Euclidean volume, then the Riemannian
curvature tensor satisfies $$|Rm|_{g(t)}(x,t)\leq t^{-1},\mbox{
whenever }{\rm{dist}}_{g(t)}(x_{0},x)<\epsilon_{0}r_{0}, \ \ \ {\rm
and } \ \ 0<t\leq (\epsilon_{0}r_{0})^{2}.$$

In particular, $|Rm|_{g(t)}(x_{0},t)\leq t^{-1}$ for all time
$t\in(0,(\epsilon_{0}r_{0})^{2}]$.
\end{thm}
\begin{proof}
We use Claim 1 and Claim 2 of Theorem 10.1 in [Pe1] and adopt a
contradiction argument. For any given small constants
$\epsilon,\delta>0$, set
$\epsilon_{0}=\epsilon,\delta_{0}=\delta$, then there is a
solution to the Ricci flow equation (1.2), say $(M,g(t))$, not
satisfying the conclusion of the theorem. After a rescaling, we
may assume that $r_{0}=1$. Denote by $\bar{M}$ the non-empty set
of pairs $(x,t)$ such that $|Rm|_{g(t)}(x,t)> t^{-1}$, then as in
Claim 1 and Claim 2 of Theorem 10.1 in [Pe1], we can choose
another space time point $(\bar{x},\bar{t})\in\bar{M}$ with
$0<\bar{t}\leq\epsilon^{2},{\rm{dist}}_{g(\bar{t})}(x_{0},\bar{x})<\frac{1}{10}$,
such that $|Rm|_{g(t)}(x,t)\leq4Q$ whenever
$$\bar{t}-\frac{1}{2n}Q^{-1}\leq t\leq\bar{t},\hspace{1cm}{\rm{dist}}_{g(\bar{t})}(\bar{x},x)\leq\frac{1}{10}(100n\epsilon)^{-1}Q^{-1/2},$$
where $Q=|Rm|_{g(\bar{t})}(\bar{x},\bar{t})$. It is remarkable
that from the proof of Claim 2 of Theorem 10.1 in [Pe1], each such
a space time point $(x,t)$ satisfies
$${\rm{dist}}_{g(t)}(x,x_{0})<{\rm{dist}}_{g(\bar{t})}(x_{0},\bar{x})+(100n\epsilon)^{-1}Q^{-1/2}<\frac{1}{10}+(100n)^{-1}<\frac{1}{2}.$$

Now choosing sequences of positive numbers
$\epsilon_{k}\rightarrow0$ and $\delta_{k}\rightarrow0$, we obtain
a sequence of solutions
$(M_{k},g_{k}(t)),t\in[0,\epsilon_{k}^{2}]$ and a sequence of
points $x_{0,k},\bar{x}_{k}\in M_{k}$ and times $\bar{t}_{k}$,
with each satisfying the assumptions of the theorem and the
properties described above. In particular, we have that
$Q_{k}=|Rm_{k}|_{g_{k}(\bar{t}_{k})}(\bar{x}_{k},\bar{t}_{k})\rightarrow\infty$.
Consider the sequence of pointed Ricci flow solutions
$$(B_{g_{k}(\bar{t}_{k})}(\bar{x}_{k},
\frac{1}{10}(100n\epsilon_{k})^{-1}Q_{k}^{-1/2}),Q_{k}g_{k}(Q_{k}^{-1}t+\bar{t}_{k}),\bar{x}_{k}),
t\in[-\frac{1}{2n},0].$$ Using Hamilton's compactness theorem for
solutions to the Ricci flow,
 we can extract a subsequence which
converge to a complete Ricci flow solution
$(M_{\infty},g_{\infty}(t),\bar{x}_{\infty}),t\in(-\frac{1}{2n},0],$
with $|Rm_{\infty}|_{g_{\infty}(0)}(\bar{x}_{\infty},0)=1$.

By assumption, the balls  $$B_{g_{k}(\bar{t}_{k})}(\bar{x}_{k},
\frac{1}{10}(100n\epsilon_{k})^{-1}Q_{k}^{-1/2})\subset
B_{g_{k}(t)}(x_{0,k},\frac{1}{2})$$ for any
$t\in[\bar{t}-\frac{1}{2n}Q^{-1},\bar{t}]$, so the scalar
curvature $R_{k}(x,t)\geq-1$ for
$t\in[\bar{t}-\frac{1}{2n}Q^{-1},\bar{t}]$ and $x\in
B_{g_{k}(\bar{t}_{k})}(\bar{x}_{k},
\frac{1}{10}(100n\epsilon_{k})^{-1}Q_{k}^{-1/2})$ and
${\rm{Vol}}_{g_{k}(t)}(B_{g_{k}(t)}(x,r))\geq(1-\delta_{k}){\rm{Vol}}(B(r))$
for any metric ball $B_{g_{k}(t)}(x,r)\subset
B_{g_{k}(\bar{t}_{k})}(\bar{x}_{k},
\frac{1}{10}(100n\epsilon_{k})^{-1}Q_{k}^{-1/2})$,
$t\in[\bar{t}-\frac{1}{2n}Q^{-1},\bar{t}]$. Passing to the limit,
we see that $g_{\infty}(t)$ has scalar curvature $R_{\infty}\geq0$
everywhere and local volume
${\rm{Vol}}_{g_{\infty}(t)}(B_{g_{\infty}(t)}(z,r))\geq
{\rm{Vol}}(B(r))$ for any balls $B_{g_{\infty}(t)}(z,r)$ at time
$t\in(-\frac{1}{2n},0]$. Then the local variation formula of
volume implies that $R_{\infty}\equiv0$ on
$M_{\infty}\times(-\frac{1}{2n},0]$, see [STW] for details. By the
evolution of the scalar curvature $\frac{\partial}{\partial
t}R_{\infty}=\triangle R_{\infty}+2|Ric_{\infty}|^{2}$, we get
that $Ric_{\infty}\equiv0$ over
$M_{\infty}\times(-\frac{1}{2n},0]$. Then the Bishop-Gromov volume
comparison theorem implies that $g_{\infty}(t)$ are flat solutions
to the Ricci flow, which contradicts the fact that
$|Rm_{\infty}|(\bar{x}_{\infty},0)=1$. This ends the proof of the
theorem.
\end{proof}

The next lemma provides a comparison of the curvature of the
normalized and unnormalized Ricci flow. By assumption, there is
$\bar{C}<\infty$ such that $|Ric|\leq\bar{C}$ everywhere along the
flow $(M,g(t))$. Note that by Lemma 3.1, there is some time
$T<\infty$ such that $2\overline{R}_{\infty}\leq
\overline{R}(g(t))\leq\frac{1}{2}\overline{R}_{\infty}<0$ whenever
$t>T$. Fix any such a time $\bar{t}>T$ and let $h(t)$ and
$\tilde{h}(\tilde{t})$ be the solutions to the normalized and
unnormalized Ricci flow with initial metric
$h(0)=\tilde{h}(0)=g(\bar{t})$ respectively, where
$\tilde{t}=\tilde{t}(t)$ is the corresponding rescaled time for
$t$. Denote by $Rm_{\bar{t}},Ric_{\bar{t}},R_{\bar{t}}$ and
$\widetilde{Rm}_{\bar{t}},\widetilde{Ric}_{\bar{t}},\tilde{R}_{\bar{t}}$
the corresponding Riemannian curvature, Ricci curvature and scalar
curvature of them, where $|Ric_{\bar{t}}|\leq\bar{C}$ since
$h(t)=g(\bar{t}+t)$. Then we have
\begin{lem}
The solution $\tilde{h}(\tilde{t})$ exists for all time
$\tilde{t}\in[0,\infty)$. Furthermore, there exist constants $C$
and $\tau$ depending on $\bar{\lambda}_{M}$ and $\bar{C}$, such
that
$$t\leq\tilde{t}\leq Ct,
|\widetilde{Rm}_{\bar{t}}|(x,\tilde{t})\leq|Rm_{\bar{t}}|(x,t)\leq
C|\widetilde{Rm}_{\bar{t}}|(x,\tilde{t}), \mbox{ whenever }
t\leq\tau.$$
\end{lem}
\begin{proof}
The solution $h(t)$ has average scalar curvature
$\overline{R}(\bar{t}+t)\leq\frac{1}{2}\overline{R}_{\infty}<0$,
so $\tilde{h}(\tilde{t})$ also has average scalar curvature
$\tilde{\overline{R}}<0$. From the evolution
$\frac{d}{d\tilde{t}}\ln
{\rm{Vol}}(\tilde{h}(\tilde{t}))=-\tilde{\overline{R}}$, the
volume ${\rm{Vol}}(\tilde{h}(\tilde{t}))$ increases strictly in
$\tilde{t}$, so to normalize it, we need to compress the space and
time. Thus $\tilde{t}\geq t$ and
$|\widetilde{Rm}_{\bar{t}}|(x,\tilde{t})\leq|Rm_{\bar{t}}|(x,t)$
for all $(x,t)$. So $\tilde{h}(\tilde{t})$ exists for all time.

The last assertion means that the scaling factor from normalized
Ricci flow to the unnormalized one is less than $C$ on the time
interval $[0,\tau]$. Consider the evolution of average scalar
curvature $\tilde{\overline{R}}(\tilde{t})$:
\begin{equation}\nonumber
\frac{d}{d\tilde{t}}\tilde{\overline{R}}=\frac{\int_{M}(2|\widetilde{Ric}_{\bar{t}}|^{2}-\tilde{R}_{\bar{t}}^{2})dv_{k}}
{Vol_{\tilde{h}^(\tilde{t})}(M)}+\tilde{\overline{R}}^{2}\leq\Lambda,
\end{equation}
for some constant $\Lambda=\Lambda(\bar{C})$, since
$|\widetilde{Ric}_{\bar{t}}|\leq|Ric_{\bar{t}}|\leq\bar{C},
|\tilde{R}_{\bar{t}}|\leq|R_{\bar{t}}|\leq\bar{C},|\tilde{\overline{R}}|\leq|\overline{R}|\leq\bar{C}$.
Note that the initial value
$\tilde{\overline{R}}(0)=\overline{R}(g(\bar{t}))\leq\frac{1}{2}\overline{R}_{\infty}$,
so there is some constant $\tilde{\tau}=\tilde{\tau}(\Lambda)$
such that
$\tilde{\overline{R}}(\tilde{t})\leq\frac{1}{4}\overline{R}_{\infty}$
for $\tilde{t}\in[0,\tilde{\tau}]$. Thus the scaling factor from
normalized Ricci flow to the unnormalized one, which equals
$\frac{\overline{R}(h(t))}{\tilde{\overline{R}}(\tilde{t})}$, is
less than $8$ on the time interval $\tilde{t}\in[0,\tilde{\tau}]$.
Now the result follows, by setting $\tau=\frac{\tilde{\tau}}{8}$
and $C=8$.
\end{proof}

The following lemma gives the estimation of the local volume along
the Ricci flow. As in [Se], the proof uses Theorem A 1.5 of [CC].
By assumption, we have a solution $(M,g(t))$ to the normalized
Ricci flow (1.3) and a sequence of times $t_{k}\rightarrow\infty$
and points $p_{k}$ such that $(M,g(t_{k}),p_{k})\stackrel
{d_{GH}}\longrightarrow (N_{\infty},g_{\infty},p_{\infty})$ with
$g(t_k) \stackrel {C^{1,\alpha}}\longrightarrow g_\infty$ on the
regular part $\mathcal{R}$ of the orbifold $N_{\infty}$. For the
space $M$ or $N_{\infty}$, let $\mathcal{R}_{\epsilon,\rho}$ be
the set of points $x$ such that $d_{GH}(B(x,r),B(r))<\epsilon r$
for any $r\leq u$, where $u\geq\rho$ is some constant depending on
$x$. Here and after, $B(r)$ denotes a ball of radius $r$ in
$4$-Euclidean space and $B(x,r)$ the metric ball of radius $r$
with center $x$ in a metric space. A weak version is
$\mathcal{WR}_{\epsilon,\rho}$, the set of points $x$ such that
there is $u\geq\rho$ with $d_{GH}(B(x,u),B(u))<\epsilon u$.

\begin{lem}
For each $q\in\mathcal{R}$, choose a sequence $q_{k}\in M$ that
converge to $q$. Then for any $\epsilon>0$, there exist
$k_{0},\eta,\rho>0$ such that
$${\rm{Vol}}(B_{g(t_{k}+t)}(q_{k}^{'},r))\geq(1-\epsilon){\rm{Vol}}(B(r)),\forall r<\rho,k_{0}<k,$$
whenever $B_{g(t_{k}+t)}(q_{k}^{'},r)\subset
B_{g(t_{k})}(q_{k},\rho)$ and $t\in[-\eta,\eta]$.
\end{lem}
\begin{proof}
By the boundedness of Ricci tensor, there is a universal constant
$\Lambda=\Lambda(\bar{C})>1$ such that
$B_{g(t)}(p,\Lambda^{-1}r)\subset B_{g(s)}(p,r)\subset
B_{g(t)}(p,\Lambda r)$ for all $t,s\in[t_{k}-1,t_{k}+1],p\in M$
and $r>0$. By Theorem A.1.5 of [CC], for fixed $\epsilon>0$, there
are $\delta=\delta(\epsilon,n),\rho=\rho(\epsilon,n)>0$ such that
$x\in\mathcal{WR}_{\delta,r}$ implies
${\rm{Vol}}(B_{g(t)}(x,r))\geq(1-\epsilon){\rm{Vol}}(B(r))$ for
each $r\leq\rho$ and $x\in M$. So by definition, it suffice to
show $q_{k}^{'}\in\mathcal{R}_{\delta,\rho}$ with respect to each
metric $g(t),t\in[t_{k}-\eta,t_{k}+\eta]$, whenever $q_{k}^{'}\in
B_{g(t_{k})}(q_{k},\Lambda\rho)$, for some constant $\eta>0$. The
constant $\rho$ may be modified by a smaller one if necessary.

Using Theorem A.1.5 of [CC] again, for fixed $\delta$ as above,
there is $\delta_{1}=\delta_{1}(\delta,n)>0$ such that
$q_{k}\in\mathcal{WR}_{\delta_{1},\frac{(\Lambda^{2}+1)\rho}{1-\delta}}$
 implies $q_{k}^{'}\in\mathcal{R}_{\delta,\rho}$ for any
$q_{k}^{'}\in B_{g(t)}(q_{k},\Lambda^{2}\rho)$. So it reduces to
show
$q_{k}\in\mathcal{WR}_{\delta_{1},\frac{(\Lambda^{2}+1)\rho}{1-\delta}}$
with respect to each time $t\in[t_{k}-\eta,t_{k}+\eta]$ for some
$\eta>0$ small enough. In fact, as showed in [Se],
$d_{GH}(B_{g(t_{k})}(q_{k},\rho_{1}),B(\rho_{1}))<\frac{1}{2}\delta_{1}\rho_{1}$
for some small number $\rho_{1}$ and all $k$ large enough. By the
boundedness of Ricci tensor again, there is a constant $\eta\leq1$
such that for each time $t\in[-\eta,\eta]$, we have
$d_{GH}(B_{g(t_{k}+t)}(q_{k},\rho_{1}),B_{g(t_{k})}(q_{k},\rho_{1}))<\frac{1}{2}\delta_{1}\rho_{1}$
for all $k$. Thus
$d_{GH}(B_{g(t_{k}+t)}(q_{k},\rho_{1}),B(\rho_{1}))<\delta_{1}\rho_{1}$
for each $t\in[-\eta,\eta]$. Now the result follows by setting
$\rho=\frac{(1-\delta)\rho_{1}}{\Lambda^{2}+1}$.
\end{proof}

Note that in the proof, the constant
$\delta_{1}=\delta_{1}(\epsilon,n)$, so the constant $\eta$
depends only on $\epsilon,n$ and $\bar{C}$. By assumption, there
is a compact exhaustion $\{K_{i}\}_{i=1}^{\infty}$ of
$\mathcal{R}$ and a sequence of smooth embeddings
$F_{i}:K_{i}\rightarrow M$ such that $F_{i}(p_{\infty})=p_{i}$ and
$F_{i}^{*}g(t_{i})$ converges to $g_{\infty}$ in the local
$C^{1,\alpha}$ sense. Following the lines described in [Se], we
can prove

\begin{lem}
   Denote by $K_{i,k}=F_{k}(K_{i})$, then for any $\epsilon>0$ and $i$, there
are $k_{0},\eta,\rho>0$ such that
$${\rm{Vol}}(B_{g(t_{k}+t)}(q_{k}^{'},r))\geq(1-\epsilon){\rm{Vol}}(B(r)),\forall q_{k}^{'}\in K_{i,k},k_{0}<k,t\in[-\eta,\eta]\mbox{ and }r<\rho.$$
\end{lem}

Now we are ready to prove the Proposition 4.1.
\begin{proof}[Proof of Proposition 4.1]
Assume that $p_{\infty}\in K_{i}$ for each $i$. Set
$\epsilon=\delta_{0}$ in the the previous lemma, where
$\delta_{0}$ is just the constant in Theorem 4.2, then for one
fixed $K_{i}$, there exist $k_{0},\eta,\rho>0$ such that
${\rm{Vol}}(B_{g(t_{k}+t)}(q,r))\geq(1-\delta_{0}){\rm{Vol}}(B(r))$
whenever $q\in K_{i,k},k_{0}<k,t\in[-\eta,\eta]$ and $r<\rho$.
Modifying $\rho$ and $\eta$ by smaller constants, we assume
$(\epsilon_{0}\rho)^{2}\leq2\eta<\tau$, where $\tau$ and
$\epsilon_{0}$ are constants in Lemma 4.3 and Theorem 4.2
respectively.

Let $h_{k}(\tilde{t})$ be the corresponding solutions to the
unnormalized Ricci flow equation with initial value
$h_{k}(0)=g(t_{k}-\eta)$, then
${\rm{Vol}}(B_{h_{k}(\tilde{t})}(q,r))\geq(1-\delta_{0}){\rm{Vol}}(B(r))$
whenever $q\in K_{i,k},r<\rho,k_{0}<k$ and $\tilde{t}$ satisfying
$t(\tilde{t})\in[0,2\eta]$, since the inequality
${\rm{Vol}}(B(q,r))\geq(1-\delta_{0}){\rm{Vol}}(B(r))$ is scale
invariant and $B_{h_{k}(\tilde{t})}\subset
B_{g(t_{k}+t(\tilde{t}))}(q,r)$ for $k$ large enough such that
$t_{k}\geq T+\eta$ for $T$ chosen as above. Denote by
$\widetilde{Rm}_{k}$ the Riemannian curvature tensor of $h_{k}$,
then by Theorem 4.2 and Lemma 4.3, we have
$$|Rm|(q,t_{k}+t)\leq C|\widetilde{Rm}_{k}|(q,\tilde{t})\leq
C(\tilde{t})^{-1}\leq C(t-t_{k}+\eta)^{-1},$$ for all $q\in
K_{i,k}$. Hence $|Rm|(q,t)$ is uniformly bounded on
$K_{i,k}\times[t_{k}-\frac{\eta}{2},t_{k}+\frac{\eta}{2}]$.

By Hamilton's compactness theorem of Ricci flow solution,
 $\{(K_{i,k},g(t_{k}+t),p_{k})\}_{k=1}^{\infty}$
converge along a subsequence to a solution to the normalized Ricci
flow
$(K_{i,\infty},g_{i,\infty}(t),p_{i,\infty}),t\in(-\frac{\eta}{2},\frac{\eta}{2})$,
in the local $C^{\infty}$ sense. When we consider the time $t=0$,
then using a diagonalization argument, a subsequence of
$\{(K_{i,k},g(t_{k}),p_{k})\}_{i,k}$ will converge in the local
$C^{\infty}$ sense to a smooth Riemannian manifold
$(K_{\infty},g_{\infty},p_{\infty})$, which is just
$(\mathcal{R},g_{\infty})$, by the uniqueness of the limit space.

For fixed $i$, there is a family of metrics
$g_{i,\infty}(t),t\in(-\frac{\eta}{2},\frac{\eta}{2})$, on
$K_{i}$. As showed in [Se], we translate the time by
$\frac{\eta}{4}$, say considering the sequence
$\{(K_{i,k},g(t_{k}+\frac{\eta}{4}+t),p_{k})\}_{k}$, and repeat
the above argument, then obtain that
$\{(K_{i,k},g(t_{k}+t),p_{k})\}_{k}\stackrel{C_{loc}^{\infty}}{\longrightarrow}(K_{i,\infty},g_{i,\infty}(t),p_{i,\infty})$
along another subsequence, on the time interval
$t\in(-\frac{\eta}{2},\frac{\eta}{4}+\frac{\eta}{2})$. The
essential point is that the estimate
$d_{GH}(B_{g(t_{k})}(q_{k},\rho_{1}),B(\rho_{1}))<\frac{1}{2}\delta_{1}\rho_{1}$
in the proof of Lemma 4.4 holds for some constant $\rho_{1}$,
simultaneously the time $t_{k}$ is replaced by
$t_{k}+\frac{\eta}{4}$, but the constant $\eta$ in Lemma 4.5 is
fixed in this procedure. Iterating this process infinite times we
obtain the convergence on $K_{i}$ for all $t\in [0,\infty)$. Then
do the same thing for each $K_{i},i=1,2,\cdots$, and after a
diagonalization argument, we get that a subsequence of
$\{(K_{i,k},g(t_{k}+t),p_{k})\}_{k}$, say
$(K_{i,k_{i}},g(t_{k_{i}}+t),p_{k_{i}})\stackrel{C_{loc}^{\infty}}{\longrightarrow}(\mathcal{R},g_{\infty}(t),p_{\infty})$
 for all $t\in
[0,\infty)$, with $g_{\infty}(0)=g_{\infty}$.

We finally show that the completion of $\mathcal{R}$ with respect
to the metric $g_{\infty}(t),$ say $\bar{\mathcal{R}}_{t}$, is
just $N_{\infty}$, for each time $t\in [0,\infty)$. Denote by
$\mathcal{S}=N_{\infty}\backslash\mathcal{R}$ the set of singular
points of $(N_{\infty},g_{\infty}(0))$, then it suffice to show
that $\bar{\mathcal{R}}_{t}=\mathcal{R}\cup\mathcal{S}$ for fixed
time $t$. Assume $\mathcal{S}=\{q_{l}\}_{l=1}^{Q}$, where
$Q\leq\beta$ for $\beta=\beta(M)$ by Lemma 3.5, and let
$\varepsilon>0$ be any small constant such that
$B_{g_{\infty}(0)}(q_{i},\varepsilon)\cap
B_{g_{\infty}(0)}(q_{j},\varepsilon)=\emptyset$ whenever $i\neq
j$. Denote by
$K_{\varepsilon}=\mathcal{R}\backslash\bigcup_{p_{l}}B_{g_{\infty}(0)}(p_{l},\varepsilon)$,
 then using $|Ric_{\infty}|\leq\bar{C}$ on $\mathcal{R}\times[0,\infty)$ and
 by the evolution of the distance function, we obtain $d_{GH}((\mathcal{R}\backslash K_{\varepsilon},g_{\infty}(t)),\mathcal{S})\leq
 e^{2\bar{C}t}\varepsilon$ and consequently
 $\bar{\mathcal{R}}_{t}=\mathcal{R}\cup\mathcal{S}$, by letting
 $\varepsilon\rightarrow0$.
\end{proof}

\section{Proofs of Theorems 1.1 and 1.2}

The main result of this section is the following

 \begin{thm}  Let $(M, \mathfrak{c})$ be a smooth oriented closed $4$-manifold
with a $\rm Spin^{c}$-structure $\mathfrak{c}$. Assume that  the
first Chern class $c_1 (\mathfrak{c})$ of $\mathfrak{c}$ is a
monopole class of $M$ satisfying that
\begin{equation}
c_{1}^{2}(\mathfrak{c})[M] \geq 2\chi(M)+3\tau(M)
>0.
\end{equation}
Let  $g(t),t\in[0,\infty)$, be   a solution to (1.3)
so that $|Ric(g(t))|\leq 3$, and
\begin{equation}
\lim\limits_{t\to \infty}\overline{\lambda}_{M}(g(t))=-
\sqrt{32\pi^{2} c_{1}^{2}(\mathfrak{c})[M]}.
\end{equation}
Then there exists an $m\in \mathbb{N}$, and  sequences  of points
$\{x_{j,k}\in M\}$, $j=1, \cdots, m$, satisfying that, by passing to
a subsequence,
$$(M, g(t_{k}+t),x_{1,k},\cdots, x_{m,k}) \stackrel{d_{GH}}\longrightarrow (
\coprod _{j=1}^m N_j , g_{\infty}, x_{1,\infty}, \cdots, ,
x_{m,\infty}),$$ $t\in [0, \infty)$,  in the $m$-pointed
Gromov-Hausdorff sense for any  $t_{k}\to \infty$, where  $(N_{j},
g_{\infty})$ $j=1,\cdots, m$ are complete K\"ahler-Einstein
orbifolds of complex  dimension 2 with at most finitely many
isolated orbifold points $\{q_{i}\}$.
 The     scalar curvature (resp. volume) of $g_{\infty}$ is
$$-\text{Vol}_{g_{0}}(M)^{-\frac{1}{2}}\sqrt{32\pi^{2}
c_{1}^{2}(\mathfrak{c})[M]}\ \ {(resp.}  \ \ \ \
V=\text{Vol}_{g_{0}}(M)=\sum_{j=1}^{m}\text{Vol}_{g_{\infty}}(N_{j})).$$
Furthermore, in the regular part of $N_{j}$, $\{g(t_{k}+t)\}$
converges to $g_{\infty}$ in  $C^\infty $-sense.
\end{thm}

Comparing with Proposition 3.2, Theorem 5.1 shows that the Einstein
orbifolds are actually  K\"ahler Einstein orbifolds  under the
additional assumptions. The key point in the proof is that the
sequence of the self-dual parts of the curvatures of the connections
on the determinant line bundles given by the irreducible solutions
in the Seiberg-Witten equations converges to a non-trivial parallel
self-dual $2$-form on every component $N_j$, which is a candidate of
the K\"ahler form.

 Let $(M, \mathfrak{c})$ and $g(t)$ be the same as in
 Thoerem  5.1, and
  let $V$,  $m$,  $t_{k}$, $x_{j,k}$, $\breve{R}(g(t))$, $g_{k}$, $g_{\infty}$,   $N_{j}$ and $F_{j,k,r}$
be the same as in Section  3. Assume that, for each $k$, $(\phi_{k},
A_{k})$ is  an irreducible solution to the Seiberg-Witten equations
(2.1). Let $|\cdot |_k$   denote the norm with respect to the metric
$g_k=g(t_{k})$. The following lemma shows that the $L^2$-norms of
the self-dual parts $F^{+}_{A_{k}}$ tends to zero.

 \begin{lem}
$$ \lim_{k\longrightarrow \infty}\int_{M}|\nabla^{k}
F^{+}_{A_{k}}|_{k}^{2}dv_{k}=0,$$ where $\nabla^{k}$ is the
connection on $\Lambda ^2T^*(M)$ induced by
 Levi-civita connection.
   \end{lem}

    \begin{proof} The Bochner formula implies that
$$0=-\frac{1}{2}\Delta_{k}
|\phi_{k}|_{k}^{2}+|\nabla^{A_{k}}\phi_{k}|_{k}^{2}+\frac{R(g_{k})}{4}|\phi_{k}|_{k}^{2}+\frac{1}{4}|\phi_{k}|_{k}^{4},$$
By taking integration we get that
\begin{equation}\int_{M}(|\nabla^{A_{k}}\phi_{k}|_{k}^{2}+\frac{R(g_{k})}{4}|\phi_{k}|_{k}^{2})dv_{k}=-\frac{1}{4}
\int_{M}|\phi_{k}|_{k}^{4}dv_{k}.\end{equation} Since
$\lambda_{M}(g_{k})$ is the lowest eigenvalue of the
 operator $-4\triangle_{k}+R(g_{k})$,
      for any $1\gg \epsilon >0$,  by definition \begin{equation}
 \lambda_{M}(g_{k})\int_{M}|\phi_{k}|_{k,\epsilon}^{2}dv_{k}\leq
  \int_{M}(4|\nabla|\phi_{k}|_{k,\epsilon}|^{2}+R(g_{k})|\phi_{k}|_{k,\epsilon}^{2})dv_{k},\end{equation}
  where
      $|\cdot|_{k,\epsilon}^{2}=|\cdot|_{k}^{2}+\epsilon^{2}$. By
      Kato's inequality (cf. (2.5)) and letting $\epsilon \rightarrow 0$, $$\lambda_{M}(g_{k})\int_{M}|\phi_{k}|_{k}^{2}dv_{k}\leq
      \int_{M}(4|\nabla^{A_{k}}\phi_{k}|_{k}^{2}+R(g_{k})|\phi_{k}|_{k}^{2})dv_{k}=-
\int_{M}|\phi_{k}|_{k}^{4}dv_{k}\leq 0. $$
 As $ \lambda_{M}(g_{k})\leq 0$, by Schwarz inequality, \begin{eqnarray*}
\overline{
\lambda}_{M}(g_{k})(\int_{M}|\phi_{k}|_{k,\epsilon}^{4}dv_{k})^{\frac{1}{2}}
  = \lambda_{M}(g_{k})\text{Vol}_{g_{k}}(M)^{\frac{1}{2}}(\int_{M}|\phi_{k}|_{k,\epsilon}^{4}dv_{k})^{\frac{1}{2}}
  & \leq & \lambda_{M}(g_{k})\int_{M}|\phi_{k}|_{k,\epsilon}^{2}dv_{k}.\end{eqnarray*}
 Therefore  $$ \overline{
\lambda}_{M}(g_{k})(\int_{M}|\phi_{k}|_{k,\epsilon}^{4}dv_{k})^{\frac{1}{2}}
   \leq
  \int_{M}(4|\nabla|\phi_{k}|_{k,\epsilon}|^{2}+R(g_{k})|\phi|_{k,\epsilon}^{2})dv_{k}.$$  Thus \begin{equation}4
\int_{M}(|\nabla^{A_{k}}\phi_{k}|_{k}^{2}-|\nabla|\phi_{k}|_{k,\epsilon}|^{2})dv_{k}\leq
- \int_{M}|\phi_{k}|_{k}^{4}dv_{k}- \overline{
\lambda}_{M}(g_{k})(\int_{M}|\phi_{k}|_{k,\epsilon}^{4}dv_{k})^{\frac{1}{2}}.\end{equation}
From (2.5), $|\nabla|\phi_{k}|_{k,\epsilon}|^{2}\leq
\frac{3}{4}|\nabla^{A_{k}}\phi_{k}|_{k}^{2}$. Hence,  by letting
$\epsilon \longrightarrow 0$, we have \begin{equation}
\int_{M}|\nabla^{A_{k}}\phi_{k}|_{k}^{2}dv_{k}\leq -(
(\int_{M}|\phi_{k}|_{k}^{4}dv_{k})^{\frac{1}{2}}+ \overline{
\lambda}_{M}(g_{k}))(\int_{M}|\phi_{k}|_{k}^{4}dv_{k})^{\frac{1}{2}}.\end{equation}
If $c_{1,k}^{+}$ denotes the self-dual part of the harmonic form
representing the first Chern class $c_{1}(\mathfrak{c})$ of
$\mathfrak{c}$, by the Seiberg-Witten equation we get that
\begin{equation}\int_{M}|\phi_{k}|_{k}^{4}dv_{k}=
8\int_{M}|F^{+}_{A_{k}}|^{2}dv_{k}\geq
32\pi^{2}[c_{1,k}^{+}]^{2}[M]\geq
32\pi^{2}c_{1}^{2}(\mathfrak{c})[M].
\end{equation}
 Note that, by the standard estimates for Seiberg-Witten
equations,  $$-\breve{R}(g_{k})\geq |\phi_{k}|_{k}^{2}$$ and, by
Theorem 1.1 in [FZ], $\sqrt{32\pi^{2}
c_{1}^{2}(\mathfrak{c})[M]}+\overline{ \lambda}_{M}(g_{k})$ is
non-positive.
 Hence\begin{equation}\end{equation} \begin{eqnarray*}\int_{M}|\nabla^{A_{k}}\phi_{k}|_{k}^{2}dv_{k}
 & \leq & -(\sqrt{32\pi^{2}
c_{1}^{2}(\mathfrak{c})[M]}+\overline{
\lambda}_{M}(g_{k}))(\int_{M}|\phi_{k}|_{k}^{4}dv_{k})^{\frac{1}{2}}\\
&\leq &  \breve{R}(g_{k})V^{\frac{1}{2}}(\sqrt{32\pi^{2}
c_{1}^{2}(\mathfrak{c})[M]}+\overline{
\lambda}_{M}(g_{k}))\longrightarrow 0,\end{eqnarray*} when
$k\longrightarrow \infty$, by (5.2) and Lemma 3.1.

     By the second one of the Seiberg-Witten
 equations again (cf. [Le2]), \begin{equation}|\nabla^{k} F^{+}_{A_{k}}|_{k}^{2}\leq
 \frac{1}{2}|\phi_{k}|_{k}^{2}|\nabla^{A_{k}} \phi_{k}|_{k}^{2},\end{equation} where $\nabla^{A_k}$ is the connection on
 $\Gamma (S_\mathfrak c)$ induced by
 the  Levi-civita connection.  Hence
 $$\int_{M}|\nabla^{k} F^{+}_{A_{k}}|_{k}^{2}dv_{k}\leq  \frac{1}{2}|\breve{R}(g(t_{k}))|\int_{M}|
\nabla^{A_{k}}\phi_{k}|_{k}^{2}dv_{k}
 \longrightarrow 0,$$ when $k\longrightarrow \infty$.
\end{proof}

Regard $F^{+}_{A_{k}}$ as self-dual 2-forms of  $g'_{k}$ on $
U_{j,r}=B_{g_{\infty}}(x_{j,\infty},r)\backslash
\bigcup_{i}B_{g_{\infty}}(q_{i,j},r^{-1})$, where
 $g'_{k}=F_{j,k,r+1}^{*}g_{k}$, and $q_{i,j} $ are the orbifold points of $N_{j}$. Since
 \begin{equation}|F^{+}_{A_{k}}|_{k}^{2}=\frac{1}{8}|\phi_{k}|_{k}^{4}\leq
 \frac{1}{8}\breve{R}(g_{k})^{2}\leq C, \end{equation} where $C$ is a constant
 independent of $k$, $F^{+}_{A_{k}}\in L^{1,2}(g'_{k})$, and $$\|F^{+}_{A_{k}}\|_{L^{1,2}(g'_{k})}\leq C',$$
 where $C'$ is a constant independent of $k$. Note that $\|\cdot
\|_{L^{1,2}(g_{\infty})}\leq 2 \|\cdot \|_{L^{1,2}(g'_{k})}$ for
$k\gg 1$ since $g'_{k}\stackrel{C^{1, \alpha}}\longrightarrow
g_{\infty}$ on $ U_{j,r}$. Thus, by passing to a subsequence,
$F^{+}_{A_{k}}\stackrel{L^{1,2}}\longrightarrow \Omega_{j}\in
L^{1,2}(g_{\infty})$, a  self-dual $2$-form with respect to
$g_{\infty}$.

\begin{lem} For any $j$, $\Omega_{j}$ is a smooth self-dual  2-form on
$U_{j,r}\backslash
\partial U_{j,r}$ such that  $\nabla^{\infty} \Omega_{j}\equiv 0$,
and $| \Omega_{j}|_{\infty}\equiv {\rm cont.}\neq 0$, where
$\nabla^{\infty}$ is the connection induced by the
 Levi-civita connection of $g_{\infty}$. Hence,  $g_{ \infty}$ is a K\"{a}hler metric with K\"{a}hler
form $\sqrt{2}\frac{\Omega_{j}}{|\Omega_{j}|}$ on $U_{j,r}$.
    \end{lem}

\begin{proof}
By Lemma 5.2
$$0 \leq \int_{U_{j,r}}|\nabla^{\infty} \Omega_{j}|_{\infty}^{2}dv_{\infty}=\lim_{k\longrightarrow \infty}
\int_{U_{j,r}}|\nabla^{\infty}
F^{+}_{A_{k}}|_{\infty}^{2}dv_{\infty} \leq 2\lim_{k\longrightarrow
\infty} \int_{M}|\nabla^{k} F^{+}_{A_{k}}|_{k}^{2}dv_{k}=0.$$  It is
easy to see that $\Omega_{j}$ is a weak solution of   the elliptic
equation $\nabla^{\infty} \Omega_{j}=0$ on $U_{j,r}$. By elliptic
equation theory, $\Omega_{j}$ is a smooth self-dual 2-form on
$U_{j,r}\backslash
\partial U_{j,r}$,  $\nabla^{\infty} \Omega_{j}\equiv 0$, and $| \Omega_{j}|_{\infty}\equiv {\rm cont.}$.

Now we claim that, for any $j$ and $r\gg 1$,
$\int_{U_{j,r}}|\Omega_{j}|_{\infty}^{2}dv_{\infty}\neq 0$. If not,
there exist $j_{s}$, $ s=1, \cdots , m_{0}$, $m_{0}\leq m$, such
that
$\int_{U_{j_{s},r}}|\Omega_{j_{s}}|_{\infty}^{2}dv_{\infty}\equiv
0$.  By Lemma 3.1,
$\overline{R}_{\infty}=\lim\limits_{k\longrightarrow
\infty}\overline{R}(g_{k})=\lim\limits_{k\longrightarrow
\infty}\breve{R}(g_{k})= \overline{\lambda}_{M}V^{-\frac{1}{2}}$,
which is the scalar curvature of $g_{\infty}$, i.e.
$\overline{R}_{\infty}=R(g_{\infty})$. Note that,  by (5.10) and
Lemma 3.7,
\begin{eqnarray*}\int_{U_{j,r}}|\Omega_{j}|_{\infty}^{2}dv_{\infty}&=&
\lim_{k\longrightarrow \infty} \int_{U_{j,r}}|
F^{+}_{A_{k}}|_{k}^{2}dv_{k}\\
&\leq &\frac{1}{8}\lim_{k\longrightarrow \infty}
\breve{R}(g_{k})^{2}\rm Vol_{g'_{k}}(U_{j,r})\\ &=&
\frac{1}{8}\overline{R}_{\infty}^{2}\rm Vol_{g_{\infty}}(U_{j,r}),
\end{eqnarray*}
    \begin{eqnarray*} \lim_{k\longrightarrow
\infty}|\int_{M}|F^{+}_{A_{k}}|_{k}^{2}dv_{k} -\sum_{j=1}^{m}
 \int_{U_{j,r}}|F^{+}_{A_{k}}|_{k}^{2}dv_{k}|& \leq &\frac{1}{8}\lim_{k\longrightarrow \infty}
\breve{R}(g_{k})^{2}\rm Vol_{g_{k}}
 (M\backslash\bigcup_{j}F_{k,j,r}(U_{j,r}))\\ &\leq  & \frac{1}{8}C\overline{R}_{\infty}^{2}\sum_{j=1}^{m}
  \rm Vol_{g_{\infty}}(N_{j}\backslash U_{j,\frac{r}{2}}
),  \end{eqnarray*} and, by Lemma 3.1,  $$\lim_{k\longrightarrow
\infty}|\int_{M}(R(g_{k})^{2}-\overline{R}_{\infty}^{2})dv_{k}|\leq
 24 \lim_{k\longrightarrow \infty}
 \int_{M}(|R(g_{k})-\overline{R}(g_{k})|+|\overline{R}_{\infty}-\overline{R}(g_{k})|)dv_{k}
=0,$$ where $C$ is a constant in-dependent of $k$.    Hence, we
obtain
\begin{eqnarray*}\overline{R}_{\infty}^{2}\sum_{j\neq j_{1},\cdots , j_{m_{0}}}\rm Vol_{g_{\infty}}(U_{j,r})&
\geq &
\sum_{j=1}^{m}\int_{U_{j,r}}8|\Omega_{j}|_{\infty}^{2}dv_{\infty}=\lim_{k\longrightarrow
\infty}\sum_{j=1}^{m}\int_{U_{j,r}}8|F^{+}_{A_{k}}|_{k}^{2}dv_{k}\\
& \geq & \lim_{k\longrightarrow
\infty}\int_{M}8|F^{+}_{A_{k}}|_{k}^{2}dv_{k}-C\overline{R}_{\infty}^{2}\sum_{j=1}^{m}
\rm Vol_{g_{\infty}}(N_{j}\backslash U_{j,\frac{r}{2}} )
\\ & \geq&  32\pi^{2}c_{1}^{2}(\mathfrak{c})[M] -C\overline{R}_{\infty}^{2}\sum_{j=1}^{m} \rm Vol_{g_{\infty}}(N_{j}\backslash U_{j,\frac{r}{2}}
).
\end{eqnarray*} The last inequality is obtained by (5.7). Thus,   by
(5.1),  $$\overline{R}_{\infty}^{2}\sum_{j\neq j_{1},\cdots ,
j_{m_{0}}}\rm Vol_{g_{\infty}}(U_{j,r})\geq
32\pi^{2}(2\chi(M)+3\tau(M))-C\overline{R}_{\infty}^{2}\sum_{j=1}^{m}
\rm Vol_{g_{\infty}}(N_{j}\backslash U_{j,\frac{r}{2}} ).$$

  By  the Chern-Gauss-Bonnet formula and the Hirzebruch signature
theorem,  $$2\chi(M)+3\tau(M)\geq
\frac{1}{4\pi^{2}}\int_{U_{k,r}}(\frac{1}{24}
R(g_{k})^{2}+2|W^{+}(g_{k})|_{k}^{2})dv_{k}
-\frac{1}{8\pi^{2}}\int_{M}|Ric\textordmasculine(g_{k})|^{2}dv_{k}.$$
By Lemma 3.1, and the fact that
$g'_{k}\stackrel{L^{2,p}}\longrightarrow g_{\infty}$  on $U_{j,r}$,
we obtain that
\begin{eqnarray*}
 \overline{R}_{\infty}^{2}\sum_{j\neq j_{1},\cdots , j_{m_{0}}}\rm Vol_{g_{\infty}}(U_{j,r})
 & \geq &\sum_{j=1}^{m} 8\int_{U_{j,r}}(\frac{\overline{R}_{\infty}^{2}}{24}+2|W^{+}(g_{\infty})|^{2}_{\infty})dv_{\infty}\\
 & & -C\overline{R}_{\infty}^{2}\sum_{j=1}^{m} \rm Vol_{g_{\infty}}(N_{j}\backslash U_{j,\frac{r}{2}}
).
\end{eqnarray*} Note that, on any $U_{j,r}$,
$j\neq j_{1}, \cdots , j_{m_{0}}$, $\nabla^{\infty} \Omega_{j}\equiv
0$, $| \Omega_{j}|_{\infty}\equiv {\rm cont.}\neq 0$, and
$\Omega_{j}$ is a self-dual 2-form. Thus $g_{\infty}$ is a
K\"{a}hler metric with K\"{a}hler form
$\sqrt{2}\frac{\Omega_{j}}{|\Omega_{j}|}$ on  $U_{j,r}$, $j\neq
j_{1}, \cdots , j_{m_{0}}$. It is well known that
$\overline{R}_{\infty}^{2}=24|W^{+}(g_{\infty})|^{2}_{\infty}$ for
K\"{a}hler metrics (cf. [Le3]).   Thus
\begin{eqnarray*}\overline{R}_{\infty}^{2}\sum_{j\neq j_{1},\cdots ,
j_{m_{0}}}\rm Vol_{ g_{\infty}}(U_{j,r})& \geq &
\overline{R}_{\infty}^{2}\sum_{j\neq j_{1},\cdots , j_{m_{0}}}\rm
Vol_{g_{\infty}}(U_{j,r})
 -C\overline{R}_{\infty}^{2}\sum_{j=1}^{m} \rm Vol_{g_{\infty}}(N_{j}\backslash U_{j,\frac{r}{2}}
)\\ & & +\sum_{j_{s}=j_{1},\cdots , j_{m_{0}}}
8\int_{U_{j,r}}(\frac{\overline{R}_{\infty}^{2}}{24}+2|W^{+}(g_{\infty})|^{2}_{\infty})dv_{\infty}\\
& \geq & \overline{R}_{\infty}^{2}\sum_{j\neq j_{1},\cdots ,
j_{m_{0}}}\rm Vol_{g_{\infty}}(U_{j,r})
 -C\overline{R}_{\infty}^{2}\sum_{j=1}^{m} \rm Vol_{g_{\infty}}(N_{j}\backslash U_{j,\frac{r}{2}}
)\\ & & +\frac{1}{3}\sum_{j_{s}=j_{1},\cdots ,
j_{m_{0}}}\overline{R}_{\infty}^{2}\rm
Vol_{g_{\infty}}(U_{j_{s},r}).
\end{eqnarray*}

  Note that, for $r\gg 1$, $$1\gg
3C\overline{R}_{\infty}^{2}\sum_{j=1}^{m} \rm
Vol_{g_{\infty}}(N_{j}\backslash U_{j,\frac{r}{2}} )\geq
\sum_{j_{s}=j_{1},\cdots , j_{m_{0}}}\overline{R}_{\infty}^{2}\rm
Vol_{g_{\infty}}(U_{j_{s},r}).$$  A contradiction. Thus, for all
$j$, $\int_{U_{j,r}}|\Omega_{j}|_{\infty}^{2}dv_{\infty}\neq 0$, and
$\nabla^{\infty} \Omega_{j}\equiv 0$, $| \Omega_{j}|_{\infty}\equiv
{\rm cont.}\neq 0$. Thus we obtain the conclusion.

\end{proof}

\vskip 3mm

\begin{proof}[Proof of Theorem 5.1]  First, assume that $diam_{g(t_{k})}(M)\longrightarrow \infty,$ when $k\longrightarrow \infty$.
  By Proposition 3.2 and Proposition 4.1,  there exists a $m\in
\mathbb{N}$, and a sequence of points $\{x_{j,k}\in M\}$, $k\in
\mathbb{N}$, $j=1, \cdots, m$, satisfying that, by passing to a
subsequence, $(M, g(t_{k}+t), x_{1,k},\cdots, x_{m,k})$, $t\in [0,
\infty)$,  converges to $\{(N_{1}, g_{\infty}, x_{1,\infty}),
\cdots, (N_{m}, g_{\infty}, x_{m,\infty})\}$ in the $m$-pointed
Gromov-Hausdorff sense, when
 $k\longrightarrow \infty$, where  $(N_{j}, g_{\infty})$ $j=1,\cdots, m$ are complete
  Einstein 4-orbifolds  with finite isolated orbifold points $\{q_{i}\}$.
 The     scalar curvature of $g_{\infty}$ is
$$\overline{R}_{\infty}=\lim\limits_{t\longrightarrow \infty}\lambda_{M}(g(t)), \ \ \ \ \ {\rm  and}  \ \ \ \
V=\text{Vol}_{g_{0}}(M)=\sum_{j=1}^{m}\text{Vol}_{g_{\infty}}(N_{j}).$$
By Lemma 5.2, $g_{\infty} $ is a K\"{a}hler-Einstein  metric in
the non-singular  part  of $\coprod\limits_{j=1}^{m}N_{j}$.
    Then by the same  arguments as  in Section 4 of [Ti], $g_{\infty}$ is actually  a
      K\"{a}hler-Einstein orbifold  metric.  Furthermore, in  the non-singular
        part  of $\coprod\limits_{j=1}^{m}N_{j}$, $\{g(t_{k}+t)\}$, $t\in [0,
\infty)$,  $C^{\infty}$-converges to
      $g_{\infty}$ by Proposition 4.1.

      If $diam_{g_{k}}(M)<C$ for a constant $C$ in-dependent of $k$,
      we can also obtain the conclusion by the similar, but much easier,
    arguments as above.

\end{proof}

\begin{thm} Let $(M, \mathfrak{c})$ be a smooth  compact closed
oriented $ 4$-manifold with a  $\rm Spin^{c}$-structure
$\mathfrak{c}$. Assume that  the first Chern class
$c_1(\mathfrak{c})$ of $\mathfrak{c}$ is a monopole class of $M$
satisfying $c_{1}^{2}(\mathfrak{c})[M]= 2\chi(M)+3\tau(M)
>0$, and $\chi(M)=3 \tau (M)$. If $M$ admits   a
solution $g(t),t\in[0,\infty)$  to (1.3) with $|R(g(t))|\leq 12$,
then
$$\lim\limits_{t\longrightarrow
\infty}\overline{\lambda}_{M}(g(t))=-
\sqrt{32\pi^{2}c_{1}^{2}(\mathfrak{c})[M]}.$$ Furthermore, if
$|Ric(g(t))|\leq 3$, the K\"ahler-Einstein metric $g_{\infty}$ in
Theorem 5.1   is a complex hyperbolic metric.
\end{thm}

\begin{proof} Let $V=Vol_{g(t)}(M)$.   By  the Chern-Gauss-Bonnet formula and the Hirzebruch signature
theorem,  \begin{equation}2\chi(M)-3\tau(M)\geq
\frac{1}{4\pi^{2}}\int_{M}(\frac{1}{24}
R(g(t))^{2}+2|W^{-}(g(t))|^{2}
-\frac{1}{2}|Ric\textordmasculine(g(t))|^{2}dv_{g(t)},
\end{equation} where $W^{-}$ is the anti-self-dual Weyl tensor.
Note that
\begin{equation}\int_{M}R(g(t))^{2}dv_{g(t)}\geq
\overline{R}(g(t))^{2}V\longrightarrow
\overline{R}_{\infty}^{2}V=\lim\limits_{t\longrightarrow
\infty}\overline{\lambda}_{M}(g(t))^{2}, \end{equation} when
$t\longrightarrow \infty$, by  Schwarz inequality and  Lemma 3.1. By
(5.11), (5.12), Lemma 3.1 and Theorem 1.1 in [FZ],
\begin{eqnarray*} 2\chi(M)-3\tau(M) & \geq & \liminf_{t\longrightarrow
\infty}\frac{1}{2\pi^{2}} \int_{M}|W^{-}(g(t))|^{2}dv_{g(t)} +
\frac{1}{96\pi^{2}}\lim\limits_{t\longrightarrow
\infty}\overline{\lambda}_{M}(g(t))^{2} \\
&\geq  & \liminf_{t\longrightarrow \infty}\frac{1}{2\pi^{2}}
\int_{M}|W^{-}(g(t))|^{2}dv_{g(t)} + \frac{1}{3}
c_{1}^{2}(\mathfrak{c})[M]\\ &=&  \liminf_{t\longrightarrow
\infty}\frac{1}{2\pi^{2}} \int_{M}|W^{-}(g(t))|^{2}dv_{g(t)} +
\frac{1}{3} (2\chi(M)+3\tau(M)).
\end{eqnarray*}
Since $\chi(M)=3\tau(M)$, we obtain
$$\lim\limits_{t\longrightarrow
\infty}\overline{\lambda}_{M}(g(t))=-\sqrt{
32\pi^{2}c_{1}^{2}(\mathfrak{c})[M]},$$ and
$$\liminf_{t\longrightarrow
\infty} \frac{1}{2\pi^{2}}\int_{M}|W^{-}(g(t))|^{2}dv_{g(t)} =0.$$

Now, assume that $|Ric(g(t))|\leq 3$. Let $t_{k}$, $N_{j}$, $g_{k}$,
and $g_{\infty}$ be the same as above. For any $j$ and compact
subset $U$ of the regular part  of $N_{j}$,
$$0\leq \int_{U}|W^{-}(g_{\infty})|_{\infty}^{2}dv_{\infty}\leq \liminf_{k\longrightarrow \infty}
\int_{M}|W^{-}(g(t_{k}))|_{k}^{2}dv_{k} = 0,$$ since
$g(t_{k})\stackrel{L^{2,p}}\longrightarrow g_{\infty}$  on $U$.
Hence $g_{\infty}$ is a  K\"{a}hler-Einstein  metric with
$W^{-}(g_{\infty})\equiv 0$. This implies that $g_{\infty}$ is a
complex hyperbolic metric (cf. [Le1]).
 The desired result follows.

\end{proof}

\begin{proof}[Proofs of Theorem 1.1 and Theorem 1.2]

 By the work of Taubes [Ta],
  if $(M, \omega)$ is a
compact symplectic manifold with $b_2^+(M)>1$, the
spin$^c$-structure induced by $\omega$ is a monopole class.
 Moreover, since  in this situation $c_1^2(\mathfrak {c}) [M]=2\chi
(M)+3\tau (M)$,  Theorem 1.1 (resp. Theorem 1.2) is an obvious
 consequence of Theorem 5.1  (resp. Theorem  5.4).
\end{proof}

\end{document}